
\newcount\secno
\newcount\prmno
\def\section#1{\vskip1truecm
               \global\def\currenvir{section}
               \global\advance\secno by1\global\prmno=0
               {\bf \number\secno. {#1}}
               \smallskip}

\def\subsection{\global\def\currenvir{subsection}
                \global\advance\prmno by1
               \smallskip  \ind{ (\number\secno.\number\prmno) }}
\def\subsec{\global\def\currenvir{subsection}
                \global\advance\prmno by1\smallskip
                { (\number\secno.\number\prmno)\ }}

\def\proclaim#1{\global\advance\prmno by 1
                {\bf #1 \the\secno.\the\prmno$.-$ }}

\long\def\th#1 \enonce#2\endth{%
   \medbreak\proclaim{#1}{\it #2}\global\def\currenvir{th}\smallskip}

\def\rem#1{\global\advance\prmno by 1
{\it #1} \the\secno.\the\prmno$.-$ }

\magnification 1250 \pretolerance=500 \tolerance=1000
\brokenpenalty=5000 \mathcode`A="7041 \mathcode`B="7042
\mathcode`C="7043 \mathcode`D="7044 \mathcode`E="7045
\mathcode`F="7046 \mathcode`G="7047 \mathcode`H="7048
\mathcode`I="7049 \mathcode`J="704A \mathcode`K="704B
\mathcode`L="704C \mathcode`M="704D \mathcode`N="704E
\mathcode`O="704F \mathcode`P="7050 \mathcode`Q="7051
\mathcode`R="7052 \mathcode`S="7053 \mathcode`T="7054
\mathcode`U="7055 \mathcode`V="7056 \mathcode`W="7057
\mathcode`X="7058 \mathcode`Y="7059 \mathcode`Z="705A
\def\spacedmath#1{\def\packedmath##1${\bgroup\mathsurround =0pt##1\egroup$}
\mathsurround#1
\everymath={\packedmath}\everydisplay={\mathsurround=0pt}}
 \spacedmath{2pt}

\def\iso{\vbox{\hbox to .8cm{\hfill{$\scriptstyle\sim$}\hfill}
\nointerlineskip\hbox to .8cm{{\hfill$\longrightarrow $\hfill}} }}
\def\sdir_#1^#2{\mathrel{\mathop{\kern0pt\oplus}\limits_{#1}^{#2}}}

\font\eightrm=cmr8 \font\sixrm=cmr6

\def\pc#1{\tenrm#1\sevenrm}
\def\tx{\kern-1.5pt -}
\def\cqfd{\kern 2truemm\unskip\penalty 500\vrule height 4pt depth 0pt width
4pt\medbreak} 
\def\no{n\up{o}\kern 2pt}
\def\ind{\par\hskip 1truecm\relax}

\font\pal=cmsy7

\def\sp#1{{\cal S}\kern-1pt\raise-1pt\hbox{\pal P}^{}_C(#1)}

\frenchspacing
\input xy
\xyoption{all}
\input amssym.def
\input amssym
\vsize = 25truecm \hsize = 16.1truecm \voffset = -.5truecm
\parindent=0cm
\baselineskip15pt \overfullrule=0pt

\vglue 2.5truecm \font\Bbb=msbm10

\centerline{\bf On the scope of validity of the norm limitation
theorem}
\par
\centerline{\bf in one-dimensional abstract local class field
theory}
\bigskip

\centerline{I.D. {\pc CHIPCHAKOV}\footnote{$^{\ast}$}{Partially
supported by Grant MM1106/2001 of the Bulgarian Foundation for
Scientific Research.}}
\par
\vskip1.truecm
\centerline{{\bf 1. Introduction}}
\par
\medskip
This paper is devoted to the study of norm groups of finite
separable extensions of primarily quasilocal fields (briefly,
PQL-fields), i.e. of $p$-quasilocal fields with respect to every
prime number $p$. It has been proved in [11] that such a field $E$
admits one-dimensional local $p$-class field theory, provided that
the $p$-component Br$(E) _{p}$ of the Brauer group Br$(E)$ is
nontrivial. This theory shows that finite abelian $p$-extensions
of $E$ are subject to exact analogues to the local reciprocity law
and the local Hasse symbol (cf. [39, Ch. 6, Theorem 8], [25, Ch.
2, 1.3] and [11, Theorems 2.1 and 2.2]), which enables one to
obtain a satisfactory description of the norm group of any abelian
finite extension of $E$. The present paper gives two sufficient
conditions for validity of the norm group equality $N(R/E) = N(R
_{\rm ab}/E)$, where $R/E$ is a finite separable extension and $R
_{\rm ab}$ is the maximal abelian subextension of $E$ in $R$. It
shows that these conditions determine to a considerable extent the
scope of validity of the norm limitation theorem for local fields
(cf. [21, Ch. 6, Theorem 8]). This is demonstrated by describing
the norm groups of finite separable extensions of Henselian
discrete valued fields whose finite extensions are strictly
primarily quasilocal, and of formally real quasilocal fields.
\par
\medskip
The basic notions of (one-dimensional) local class field theory,
used in the sequel, are defined in Section 2. For each field $E$,
we denote by $P(E)$ the set of those prime numbers $p$, for which
there exists at least one cyclic extension of $E$ of degree $p$.
Clearly, a prime number $p$ lies in $P(E)$ if and only if $E$ does
not equal its maximal $p$-extension $E (p)$ in a separable closure
$E _{\rm sep}$ of $E$. Let us note that a field $E$ is said to be
$p$-quasilocal, if it satisfies some of the following two
conditions: (i) Br$(E) _{p} = \{0\}$ or $p \not\in P(E)$; (ii)
cyclic extensions of $E$ of degree $p$ embed as $E$-subalgebras in
each central division $E$-algebra of Schur index $p$. When this
occurs, we say that $E$ is strictly $p$-quasilocal, provided that
Br$(E) _{p} \neq \{0\}$ or $p \not\in P(E)$. The field $E$ is
called strictly primarily quasilocal, if it is strictly
$p$-quasilocal, for every prime $p$; it is said to be quasilocal,
if its finite extensions are PQL-fields. It has been proved in
[11] that strictly PQL-fields admit local class field theory. As
to the converse, it holds, in each of the following special cases:
(i) $E$ contains a primitive $p$-th root of unity, for each prime
$p \ge 5$ not equal to char$(E)$; (ii) $E$ is an algebraic
extension of a global field $E _{0}$. It should also be noted that
all presently known fields with local class field theory are
strictly PQL (see Proposition 2.6, Corollary 2.7 and the remarks
between them, for more details).
\par
\medskip
The description of norm groups of finite extensions of strictly
PQL-fields is a major objective of local class field theory. The
discussion of this problem in the classical case of a local field
$E$ usually begins with the observation that then the norm group
$N(R/E)$ of any finite extension $R$ of $E$ in $E _{\rm sep}$ is
closed of finite index in the multiplicative group $E ^{\ast }$.
Hence, by the existence theorem (cf. [21, Ch. 6, Theorem 8] or
[18, (6.2)]), $N(R/E) = N(R _{1}/E)$, for some finite abelian
extension $R _{1}$ of $E$ in $E _{\rm sep}$, uniquely determined
by $R$. As noted above, by the norm limitation theorem, we have $R
_{1} = R _{\rm ab}$. The theorem has been generalized by Moriya
[32] in the case where $R$ is separable over $E$ and $E$ possesses
a Henselian discrete valuation whose residue field $\widehat E$ is
quasifinite, i.e. perfect with an absolute Galois group $G
_{\widehat E}$ isomorphic to the profinite completion of the
additive group $\hbox{\Bbb Z}$ of integer numbers, as a part of
the development of local class field theory in this direction (see
[40; 49] and [18, Ch. V] as well). The group $N(R/E)$ also admits
an accomplished description in case $E$ is an algebraic strictly
PQL-extension  of a global field $E _{0}$. The description has
been obtained as a result of the fact that then $E$ possesses a
certain characteristic system $\{v(p): p \in P(E)\}$ of nontrivial
absolute values (see Proposition 2.8). The properties of this
system have been used in [13] for proving the validity of the
Hasse norm principle for $R/E$, and also, for showing that $N(R/E)
= N(\Phi (R)/E)$, for some finite abelian extension $\Phi (R)$ of
$E$ in $E _{\rm sep}$, uniquely determined by the local behaviour
of $R/E$ at $v(p)$ when $p$ runs through the set of elements of
$P(E)$ dividing the degree $[M:E]$ of the normal closure $M$ of
$R$ in $E _{\rm sep}$ over $E$. These results motivate one to try
to reduce the study of $N(R/E)$, for an arbitrary strictly
PQL-field $E$, to the special case in which $R$ is abelian over
$E$, and then to obtain information about $N(R/E)$, by applying
the generalization of the local reciprocity law in [11].
\par
\medskip
The purpose of this paper is to shed light on the possibility for such a
reduction by proving the following:
\par
\medskip
{\bf Theorem 1.1.} {\it Let $E$ be a field, $M/E$ a finite Galois
extension, and $R$ an intermediate field of $M/E$. Then $N(R/E) =
N(R _{\rm ab}/E)$ in each of the following special cases:}
\par
(i) $E$ {\it is primarily quasilocal and the Galois group $G(M/E)$ is
nilpotent;}
\par
(ii) {\it $E$ is a quasilocal field, such that the natural homomorphism of
Br$(E)$ into Br$(F)$ is surjective, for every finite extension of $E$.}
\par
\medskip
{\bf Theorem 1.2.} {\it For each nonnilpotent finite group $G$,
there exists a strictly} PQL{\it -field $E(G)$ and a Galois
extension $M(G)$ of $E(G)$, for which the following is true:}
\par
(i) $E(G)$ {\it is an algebraic extension of the field $\hbox{\Bbb Q}$ of
rational numbers;}
\par
(ii) {\it The Galois group $G(M(G)/E(G))$ is isomorphic to $G$, and
$N(M(G)/E(G))$ is a proper subgroup of $N(M(G) _{\rm ab}/E(G))$.}
\par
\medskip
{\bf Theorem 1.3.} {\it There exists a field $E$,  for which the
following assertions hold true:}
\par
(i) {\it finite extensions of $E$ are strictly} PQL{\it -fields;}
\par
(ii) {\it the absolute Galois group $G _{K}$ is not pronilpotent;}
\par
(iii) {\it every finite extension $R$ of $K$ is subject to the following
alternative:}
\par
($\alpha $) $R$ {\it is an intermediate field of a finite Galois extension
$M(R)/K$ with a nilpotent Galois group;}
\par
($\beta $) $N(R/K)$ {\it does not equal the norm group of any abelian finite
extension of $K$.}
\par
\medskip
Throughout the paper, algebras are understood to be associative
with a unit, simple algebras are supposed to be finite-dimensional
over their centres, Brauer groups of fields are viewed as
additively presented, homomorphisms of profinite groups are
assumed to be continuous, and Galois groups are regarded as
profinite with respect to the Krull topology. For each algebra
$A$, we consider only subalgebras of $A$ containing its unit, and
denote by $A ^{\ast }$ the multiplicative group of $A$. As usual,
a field $E$ is said to be formally real, if $-1$ cannot be
presented as a finite sum of squares of elements of $E$; we say
that $E$ is a nonreal field, otherwise. The field $E$ is called
Pythagorean, if it is formally real, and the set $E ^{\ast 2} =
\{a ^{2}: a \in E ^{\ast }\}$ is additively closed (see also [14,
Satz 1]). Our basic terminology and notation concerning valuation
theory, simple algebras, Brauer groups and abstract abelian groups
is standard (such as can be found, for example, in [15; 22; 47;
34; 29] and [19], as well as those concerning profinite groups,
Galois cohomology, field extensions, Galois theory and formally
real fields (see, for example, [42; 24] and [26]). We refer the
reader to [47, Sect. 1] and [9, Sect. 2], for the definitions of a
symbol algebra and of a symbol $p$-algebra (see also [43, Ch. XIV,
Sects. 2 and 5]).
\par
\medskip
The paper is organized as follows: Section 2 includes
preliminaries needed for the further discussion. The validity of
Theorem 1.1 in cases (i) and (ii) is established in Sections 3 and
5, respectively. In Section 4 we prove Theorem 1.2, on the basis
of the study of algebraic strictly PQL-extensions of global
fields, carried out in [13]. Section 6 contains a description of
the norm groups of Henselian discrete valued fields whose finite
extensions are strictly PQL, and a characterization of the fields
from this class with the properties required by Theorem 1.3. In
Section 7 we describe along the same lines the norm groups of formally real
quasilocal fields.
\par
\vskip1.2truecm \centerline {\bf 2. Preliminaries}
\par
\medskip
Let $E$ be a field, $E _{\rm sep}$ a separable closure of $E$,
Nr$(E)$ the set of norm groups of finite extensions of $E$, and
$\Omega (E)$ the set of finite abelian extensions of $E$ in $E
_{\rm sep}$. We say that $E$ admits (one-dimensional) local class
field theory, if the mapping $\pi $ of $\Omega (E)$ into Nr$(E)$
defined by the rule $\pi (F) = N(F/E): F \in \Omega (E)$, is
injective and satisfies the following two conditions, for each
pair $(M _{1}, M _{2}) \in \Omega (E) \times \Omega (E)$:
\par
The norm group of the compositum $M _{1}M _{2}$ is equal to the
intersection $N(M _{1}/E)$ $\cap N(M _{2}/E)$ and $N((M _{1} \cap
M _{2})/E)$ equals the inner group product $N(M _{1}/E)N(M
_{2}/E)$.
\par
We call $E$ a field with one-dimensional local $p$-class field
theory, for some prime number $p$, if the restriction of $\pi $ on
the set of abelian finite $p$-extensions of $E$ in $E _{\rm sep}$
has the same properties. Our approach to the study of fields with
local class field theory is based on the possibility of reducing
the description of norm groups of finite Galois extensions with
nilpotent Galois groups to the special case of $p$-extensions.
This possibility can be seen from the following two lemmas.
\par
\medskip
{\bf Lemma 2.1.} {\it Let $E$ be a field and $L$ an extension of
$E$ presentable as a compositum of extensions $L _{1}$ and $L
_{2}$ of $E$ of relatively prime degrees. Then $N(L/E) = N(L
_{1}/E) \cap N(L _{2}/E)$, $N(L _{1}/E) = E ^{\ast } \cap N(L/L
_{2})$, and there is a group isomorphism $E ^{\ast }/N(L/E) \cong
(E ^{\ast }/N(L _{1}/E)) \times (E ^{\ast }/N(L _{2}/E))$.}
\par
\medskip
{\it Proof.} The inclusion $N(L/E) \subseteq N(L _{i}/E)$: $i = 1,
2$ follows at once from the transitivity of norm mappings in
towers of field extensions of finite degrees (cf. [26, Ch. VIII,
Sect. 5]). Conversely, let $c \in N(L _{i}/E)$ and $[L _{i}:E] = m
_{i}$: $i = 1, 2$. As g.c.d.$(m _{1}, m _{2}) = 1$, this implies
consecutively that $[L:L _{2}] = m _{1}$, $c ^{m _{1}} \in
N(L/E)$, $[L:L _{1}] = m _{2}$, $c ^{m _{2}} \in N(L/E)$ and $c
\in N(L/E)$, and so proves the equality $N(L/E)$ $= N(L _{1}/E)
\cap N(L _{2}/E)$. Since $E ^{\ast }/N(L _{i}/E)$ is a group of
exponent dividing $m _{i}$: $i = 1, 2$, it is also clear that $N(L
_{1}/E)N(L _{2}/E) = E ^{\ast }$. These observations prove the
concluding assertion of the lemma. Our argument also shows that $N
_{L _{2}} ^{L} (\lambda _{1}) = N _{E} ^{L _{1}} (\lambda _{1})$:
$\lambda _{1} \in L _{1}$, whence $N(L _{1}/E) \subseteq E ^{\ast
} \cap N(L/L _{2})$. Considering now an element
\par \noindent
$s \in E ^{\ast } \cap N(L/L _{2})$, one obtains that $s ^{m _{2}}
\in N(L/E)$. This means that $s ^{m _{2}} \in N(L _{1}/E)$, and
since $s ^{m _{1}} \in N(L _{1}/E)$, finally yields $s \in N(L
_{1}/E)$, which completes the proof of Lemma 2.1.
\par
\medskip
{\bf Lemma 2.2.} {\it Let $E$ be a field, $M$ a finite Galois
extension of $E$ with a nilpotent Galois group $G(M/E)$, $R$ an
intermediate field of $M/E$ not equal to $E$, $P(R/E)$ the set of
prime numbers dividing $[R:E]$, $M _{p}$ the maximal $p$-extension
of $E$ in $M$, and $R _{p}$ the intersection $R \cap M _{p}$,
$\forall p \in P(R/E)$. Then the following is true:}
\par
(i) $R$ {\it is equal to the compositum of the fields $R _{p}: p \in
P(R/E)$, and
\par \noindent
$[R:E] = \prod _{p \in P(R/E)} [R _{p}:E]$;}
\par
(ii) {\it The group $N(R/E)$ equals the intersection $\cap _{p \in
P(R/E)} N(R _{p}/E)$ and $E ^{\ast }/N(R/E)$ is isomorphic to the
direct product of the groups $E ^{\ast }/N(R _{p}/E): p \in
P(R/E)$.}
\par
\medskip
{\it Proof.} Statement (i) follows from Galois theory and the
Burnside-Wielandt characterization of nilpotent finite groups (cf.
[23, Ch. 6, Sect. 2]). Proceeding by induction on the number $s$
of the elements of $P(R/E)$, and taking into account that if $s
\ge 2$, then $R _{p}$ and the compositum $R _{p} ^{\prime }$ of
the fields $R _{p'}$: $p ^{\prime } \in (P(R/E) \setminus \{p\})$,
are of relatively prime degrees over $E$, one deduces Lemma 2.2
(ii) from Lemma 2.1.
\par
\medskip
It is clear from Lemma 2.2 that a field $E$ admits local class
field theory if and only if it is a field with local $p$-class
field theory, for every $p \in P(E)$. The following lemma, proved
in [12, Sect. 4], shows that the group Br$(E) _{p}$ is necessarily
nontrivial, if $E$ admits local $p$-class field theory, for a
given $p \in P(E)$.
\par
\medskip
{\bf Lemma 2.3.} {\it Let $E$ be a field, such that} Br$(E) _{p} =
\{0\}${\it , for some prime number p. Then} Br$(E _{1}) _{p} =
\{0\}$ {\it and $N(E _{1}/E) = E ^{\ast }$, for every finite
extension $E _{1}$ of $E$ in $E (p)$.}
\par
\medskip
For convenience of the reader, we prove the following lemma, which
plays an essential role not only in the proof of Theorem 1.1 (i),
but also in the study of the norm groups of the quasilocal fields
considered in Sections 6 and 7.
\par
\medskip
{\bf Lemma 2.4.} {\it For a field $E$ and a prime number $p$, the
following conditions are equivalent:}
\par
(i) Br$(E ^{\prime }) _{p} = \{0\}${\it , for every algebraic extension $E
^{\prime }$ of $E$;}
\par
(ii) {\it The exponent of the group $E _{1} ^{\ast }/N(E _{2}/E
_{1})$ is not divisible by $p$, for any pair $(E _{1}, E _{2})$ of
finite extensions of $E$ in $E _{\rm sep}$, such that $E _{1}
\subseteq E _{2}$.}
\par
\medskip
{\it Proof.} (i)$\to $(ii): Denote by $E _{2} ^{\prime }$ the
normal closure of $E _{2}$ in $E _{\rm sep}$ over $E _{1}$, and by
$E _{1} ^{\prime }$ the intermediate field of $E _{2} ^{\prime }/E
_{1}$ corresponding by Galois theory to some Sylow $p$-subgroup of
$G(E _{2} ^{\prime }/E _{1})$. By Lemma 2.3, then we have $N(E
_{2} ^{\prime }/E _{1} ^{\prime }) = E _{1} ^{\prime \ast }$,
which implies that $N(E _{2} ^{\prime }/E _{1}) = N(E _{1}
^{\prime }/E _{1})$. This means that $E _{1} ^{\ast }/N(E _{2}
^{\prime }/E _{1})$ is of exponent dividing $[E _{1} ^{\prime }:E
_{1}]$, and since $N(E _{2} ^{\prime }/E _{1}) \subseteq N(E
_{2}/E _{1})$, $E _{1} ^{\ast }/N(E _{2}/E _{1})$ is a homomorphic
image of $E _{1} ^{\ast }/N(E _{2} ^{\prime }/E _{1})$, so its
exponent also divides $[E _{1} ^{\prime }:E _{1}]$. The obtained
result proves the implication (i)$\to $(ii).
\par
(ii)$\to $(i): Suppose for a moment that there is an algebraic
extension $F$ of $E$, such that Br$(F) _{p} \neq \{0\}$. This
implies the existence of a finite separable extension $F _{0}$ of
$E$ in $F$, possessing a noncommutative central division algebra
$\Delta _{0}$ of $p$-primary dimension (cf. [42, Ch. II, 2.3] and
[12, (1.3)]). In view of [29, Sect. 4, Theorem 2], $\Delta _{0}$
can be chosen so that ind$(\Delta _{0}) = p$. Hence, by the
Proposition in [34, Sect. 15.2], there exists a finite separable
extension $F _{1}$ of $F _{0}$ of degree not divisible by $p$, for
which $\Delta _{0} \otimes _{F _{0}} F _{1}$ is a cyclic division
$F _{1}$-algebra. This, however, leads to the conclusion that
there is a cyclic extension $F _{2}$ of $F _{1}$ of degree $p$,
embeddable in $\Delta _{0} \otimes _{F _{0}} F _{1}$ as an $F
_{1}$-subalgebra, and therefore, having the property that $N(F
_{2}/F _{1}) \neq F _{1} ^{\ast }$ (cf. [34, Sect. 15.1,
Proposition b]). The obtained contradiction completes the proof of
Lemma 2.4.
\par
\medskip
Let now $\Phi $ be a field and $\Phi _{p}$ the extension of $\Phi
$ in $\Phi _{\rm sep}$ generated by a primitive $p$-th root of
unity $\varepsilon _{p}$, for some prime number $p$. It is
well-known that then $\Phi _{p}/\Phi $ is a cyclic extension of
degree $[\Phi _{p}:\Phi ] := m$ dividing $p - 1$ (cf. [26, Ch.
VIII, Sect. 3]). Denote by $\varphi $ some $\Phi $-automorphism of
$\Phi _{p}$ of order $m$, fix an integer $s$ so that $\varphi
(\varepsilon _{p}) = \varepsilon _{p} ^{s}$, and put $V  _{i} =
\{\alpha _{i} \in \Phi _{p} ^{\ast }: \varphi (\alpha _{i})\alpha
_{i} ^{-s ^{i}} \in \Phi _{p} ^{\ast p}\}$, and $\overline V _{i}
= V _{i}/\Phi _{p} ^{\ast p}$: $i = 0,..., m - 1$. Clearly, the
quotient group $\Phi _{p} ^{\ast }/\Phi _{p} ^{\ast p} :=
\overline \Phi _{p} ^{\ast p}$ can be viewed as a vector space
over the field $\hbox{\Bbb F} _{p}$ with $p$ elements. Considering
the linear operator $\bar \varphi $ of $\overline \Phi _{p}$,
defined by the rule $\bar \varphi (\alpha \Phi _{p} ^{\ast p}) =
\varphi (\alpha )\Phi _{p} ^{\ast p}$: $\alpha \in \Phi _{p}
^{\ast }$, and taking into account that the subspace of $\overline
\Phi _{p}$, spanned by its elements $\bar \varphi ^{i} (\bar
\alpha )$: $i = 0,..., m - 1$, is finite-dimensional and $\bar
\varphi $-invariant, $\forall \bar \alpha \in \overline \Phi
_{p}$, one obtains from Maschke's theorem the following statement:
\par
\medskip
(2.1) The sum of the subspaces $\overline V _{i}$: $i = 0,..., m -
1$ is direct and equal to $\overline \Phi _{p}$.
\par
\medskip
Let $L$ be an extension of $\Phi _{p}$ in $\Phi _{\rm sep}$,
obtained by adjoining a $p$-th root $\eta _{p}$ of an element
$\beta \in (\Phi _{p} ^{\ast } \setminus \Phi _{p} ^{\ast p})$. It
is easily deduced from Kummer's theory that $[L:\Phi ] = pm$ and
the following assertions hold true:
\par
\medskip
(2.2) $L/\Phi $ is a Galois extension if and only if $\beta \in V
_{j}$, for some index $j$. When this occurs, every $\Phi _{p}$-automorphism
$\psi $ of $L$ of order $p$ satisfies the equality $\varphi ^{\prime }\psi
\varphi ^{\prime -1} = \psi ^{s'}$, where $s ^{\prime } = s ^{1-j}$ and
$\varphi ^{\prime }$ is an arbitrary automorphism of $L$ extending $\varphi $.
Moreover, the following assertions hold true:
\par
(i) $L/\Phi $ is cyclic if and only if $\beta \in V _{1}$ (Albert, see
[1, Ch. IX, Theorem 6]);
\par
(ii) $L$ is a root field over $\Phi $ of the binomial $X ^{p} - a$, for some
$a \in \Phi ^{\ast
}$, if and only if $\beta \in V _{0}$, i.e. $s ^{\prime } = s$; when this
occurs, one can take as $a$ the norm $N _{\Phi } ^{\Phi _{p}} (\beta )$.
\par
\medskip
Statements (2.1), (2.2) and the following observations will be
used for proving Theorem 1.1 (ii), as well as for describing the
norm groups of the quasilocal fields considered in Sections 6 and
7:
\par
\medskip
(2.3) For a symbol $\Phi _{p}$-algebra $A _{\varepsilon _{p}}
(\alpha , \beta ; \Phi _{p})$ (of dimension $p ^{2}$), where
$\alpha \in \Phi _{p} ^{\ast }$ and $\beta \in V _{j} \setminus
\Phi _{p} ^{\ast p}$, the following conditions are equivalent:
\par
(i) $A _{\varepsilon _{p}} (\alpha , \beta ; \Phi _{p})$ is $\Phi
_{p}$-isomorphic to $D \otimes _{\Phi } \Phi _{p}$, for some central simple
$\Phi $-algebra $D$;
\par
(ii) If $\alpha = \prod _{i=0} ^{m-1} \alpha _{i}$ and $\alpha _{i} \in V
_{i}$, for each index $i$, then $A _{\varepsilon _{p}} (\alpha , \beta ; \Phi
_{p})$ is isomorphic to the symbol $\Phi _{p}$-algebra $A _{\varepsilon _{p}}
(\alpha _{j'}, \beta ; \Phi _{p})$, where $j ^{\prime }$ is determined so that
$m$ divides $j ^{\prime } + j - 1$;
\par
(iii) With notations being as in (ii), $\alpha _{i} \in N(L/\Phi _{p})$, for
every index $i \neq j ^{\prime }$.
\par
\medskip
The main result of [12] used in the present paper can be stated
as follows:
\par
\medskip
{\bf Proposition 2.5.} {\it Let $E$ be a strictly $p$-quasilocal
field, such that $E (p) \neq E$, for some prime number $p$. Assume
also that $R$ is an extension of $E$ in $E (p)$, and $D$ is a
central division $E$-algebra of $p$-primary dimension. Then $R$ is a
$p$-quasilocal field and the following statements are true:}
\par
(i) {\it $D$ is a cyclic $E$-algebra and} ind$(D) = $ exp$(D)${\it ;}
\par
(ii) Br$(R) _{p}$ {\it is a divisible group unless $p = 2$, $R =
E$ and $E$ is a formally real field; in the noted exceptional
case,} Br$(R) _{2}$ {\it is of order $2$;}
\par
(iii) {\it The natural homomorphism of} Br$(E)$ {\it into} Br$(R)$
{\it maps} Br$(E) _{p}$ {\it surjectively on} Br$(R) _{p}${\it ;
in particular, every $E$-automorphism of the field $R$ is
extendable to a ring automorphism on each central division
$R$-algebra of $p$-primary dimension;}
\par
(iv) {\it $R$ is embeddable in $D$ as an $E$-subalgebra if and only if the
degree $[R:E]$ divides} ind$(D)${\it ; $R$ is a splitting field of $D$ if and
only if $[R:E]$ is infinite or divisible by} ind$(D)${\it ;}
\par
(v) Br$(R) _{p} = \{0\}${\it , provided that $R/E$ is an infinite extension.}
\par
\medskip
The place of strictly PQL-fields in one-dimensional local class
field theory is clarified by the following result of [11]:
\par
\medskip
{\bf Proposition 2.6.} {\it Strictly} PQL{\it -fields admit local
class field theory.
\par
Conversely, a field $E$ admitting local class field theory and satisfying the
condition} Br$(E) \neq \{0\}$ {\it is strictly} PQL{\it , provided that every
central division $E$-algebra of prime exponent $p$ is similar to a tensor
product of cyclic division $E$-algebras of Schur index $p$.}
\par
\medskip
{\bf Remark 2.7.} The question of whether central division
algebras of exponent $p$ over an arbitrary field $E$ are
necessarily similar to such tensor products is a major open
problem in the theory of central simple algebras. It is well-known
that its answer is affirmative in each of the following two
special cases: (i) if $E$ contains a primitive $p$-th root of
unity or $p = $ char$(E)$ (cf. [30, (16.1)] and [2, Ch. VII,
Theorem 30]); (ii) if $E$ is an algebraic extension of a global
field (cf. [3, Ch. 10, Corollary to Theorem 5]). Also, it has been
proved in [11] that finite extensions of a field $E$ admit local
class field theory if and only if these extensions are strictly
PQL-fields.
\par
\medskip
Our next result characterizes fields with local class field theory and with
proper maximal abelian extensions, in the class of algebraic extensions of
global fields:
\par
\medskip
{\bf Proposition 2.8.} {\it Let $E _{0}$ be a global field,
$\overline E _{0}$ an algebraic closure of $E _{0}$, and $E$ an
extension of $E _{0}$ in $\overline E _{0}$, such that $P(E) \neq
\phi $. Then the following conditions are equivalent:}
\par
(i) $E$ {\it admits local class field theory;}
\par
(ii) {\it For each $p \in P(E)$,} Br$(E) _{p}$ {\it is nontrivial
and there exists a nontrivial valuation $v(p)$ of $E$, such that
the tensor product $E (p) \otimes _{E} E _{v(p)}$ is a field,
where $E _{v(p)}$ is the completion of $E$ with respect to the
topology induced by $v(p)$.
\par
When these conditions are in force, the valuation $v(p)$ is
uniquely determined, up-to an equivalence, the natural
homomorphism of} Br$(E)$ {\it into} Br$(E _{v(p)})$ {\it maps}
Br$(E) _{p}$ {\it bijectively on} Br$(E _{v(p)}) _{p}${\it , and
$E (p) \otimes _{E} E _{v(p)}$ is isomorphic as an $E
_{v(p)}$-algebra to the maximal $p$-extension $E _{v(p)} (p)$ of
$E _{v(p)}$, for every $p \in P(E)$.}
\par
\medskip
{\bf Definition 2.9.} Let $E$ be a strictly PQL-extension of a
global field $E _{0}$, such that $P(E) \neq \phi $. By a
characteristic system of $E$, we mean a system $V(E) = \{v(p):$ $p
\in P(E)\}$ of absolute values of $E$, determined in accordance
with Proposition 2.8 (ii).
\par
\medskip
Note finally that if $E$ is an algebraic extension of a global
field $E _{0}$, and $R$ is a finite extension of $E$ in $E _{\rm
sep}$, then the group $N _{\rm loc} (R/E)$ of local norms of $R/E$
consists of the elements of $E ^{\ast }$ lying in the norm groups
$N(R _{v'}/E _{v})$, whenever $v$ is a nontrivial absolute value
of $E$, and $v ^{\prime }$ is a prolongation of $v$ on $R$. It has
been proved in [13] that if $E$ is a strictly PQL-field with $P(E)
\neq \phi $, and $R ^{\prime }$ is the normal closure of $R$ in $E
_{\rm sep}$ over $E$, then $N _{\rm loc} (R/E) \subseteq N(R/E)$
and both groups are fully determined by the local behaviour of
$R/E$ at the subset of $V(E)$, indexed by the elements of $P(E)$
dividing the degree $[R ^{\prime }:E]$, where $R ^{\prime }$ is
the normal closure of $R$ in $E _{\rm sep}$ over $E$. In this
paper, we will need this result only in the special case where $R
= R ^{\prime }$, i.e. $R/E$ is a Galois extension.
\par
\medskip
{\bf Proposition 2.10.} {\it Assume that $E _{0}$ is a global
field, $E$ is an algebraic strictly} PQL{\it -extension of $E
_{0}$ with $P(E) \neq \phi $, and $V(E) = \{v(p): p \in P(E)\}$ is
a characteristic system of $E$. Also, let $M$ be a finite Galois
extension of $E$, and $P(M/E)$ the set of prime numbers dividing
$[M:E]$. Then there exists a finite abelian extension $\widetilde
M$ of $E$ satisfying the following conditions:}
\par
(i) {\it The norm groups $N(\widetilde M/E)$, $N(M/E)$ and $N
_{\rm loc} (M/E)$ are equal;}
\par
(ii) {\it The degree $[\widetilde M:E]$ divides $[M:E]$; in particular,
$\widetilde M = E$, provided that $E (p) = E$, $\forall p \in P(M/E)$;}
\par
(iii) {\it For each prime number $p$ dividing $[\widetilde M:E]$, the
maximal
$p$-extension $\widetilde M _{p}$ of $E$ in $\widetilde M$ has the property
that $\widetilde M _{p} \otimes _{E} E _{v(p)}$ is $E _{v(p)}$-isomorphic to
the maximal abelian $p$-extension of $E _{v(p)}$ in the completion $M
_{v(p)'}$, where $v(p) ^{\prime }$ is an absolute value of $M$ extending
$v(p)$.
\par
The field $\widetilde M$ is uniquely determined by $M$, up-to an
$E$-isomorphism.}
\par
\medskip
It is worth mentioning that if $M$ is an algebraic extension of a
global field $E _{0}$, and $E$ is a subfield of $M$, such that $E
_{0} \subseteq E$ and $M/E$ is a finite Galois extension, then
$N(M/E) \subseteq N _{\rm loc} (M/E)$. Identifying $M$ with its
$E$-isomorphic copy in $E _{v,{\rm sep}}$, for a fixed nontrivial
absolute value $v$ of $E$, one deduces this from the fact that the
Galois groups $G(M _{v'}/E _{v})$ and $G(M/(M \cap E _{v}))$ are
canonically isomorphic, $N(M/E) \subseteq N(M/(M \cap E _{v}))$,
and $N _{(M \cap E _{v})} ^{M} (\mu ) = N _{E _{v}} ^{M _{v'}}
(\mu )$, in case $\mu \in M ^{\ast }$ and $v ^{\prime }$ is a
prolongation of $v$ on $M$. Moreover, it follows from Tate's
description of $N _{\rm loc} (M/E)/N(M/E)$ [6, Ch. VII, Sect.
11.4] (see also [35, Sect. 6.3]), in the special case where $E
_{0}$ is an algebraic number field and $E = E _{0}$, that $M/E$
can be chosen so that $N(M/E) \neq N _{\rm loc} (M/E)$.
\par
\vskip0.6truecm \centerline {\bf 3. Norm groups of intermediate
fields of finite}
\par
\centerline {\bf normal extensions with nilpotent Galois groups}
\par
\medskip
The purpose of this Section is to prove Theorem 1.1 (i). Clearly,
our assertion can be deduced from Galois theory, Lemma 2.2 and the
following result:
\par
\medskip
{\bf Theorem 3.1.} {\it Let $E$ be a $p$-quasilocal field, $M/E$ a
finite $p$-extension of $E$, and $R$ an intermediate field of
$M/E$. Then $N(R/E) = N(R _{\rm ab}/E)$.}
\par
\medskip
{\it Proof.} In view of Lemma 2.3 and Proposition 2.5, one may
consider only the special case in which $E (p) \neq E$ and Br$(E)
_{p}$ is an infinite group. Suppose first that $R$ is a Galois
extension of $E$. It follows from Galois theory that then the
maximal abelian extension in $R$ of any normal extension of $E$ in
$R$ is itself normal over $E$ and contains $R _{\rm ab}$ as a
subfield. Since the intermediate fields of $R/E$ are
$p$-quasilocal fields, these observations show that it is
sufficient to prove the equality $N(R/E) = N(R _{\rm ab}/E)$, under
the hypothesis that $G(R/E)$ is a Miller-Moreno group, i.e. a
nonabelian group with abelian proper subgroups. For convenience of
the reader, we begin the consideration of this case with the
following elementary lemma:
\par
\medskip
{\bf Lemma 3.2.} {\it Assume that $P$ is a Miller-Moreno
$p$-group. Then the following is true:}
\par
(i) {\it The commutator subgroup $[P, P]$ of $P$ is of order $p$, the centre
$Z(P)$ of $P$ equals the Frattini subgroup $\Phi (P)$, and the group $P/\Phi
(P)$ is elementary abelian of order $p ^{2}$;}
\par
(ii) {\it A subgroup $H$ of $P$ is normal in $P$ if and only if $[P,P]
\subseteq
H$ or $H \subseteq \Phi (P)$;}
\par
(iii) {\it The quotient group of $P$ by its normal subgroup $H _{0}$ is
cyclic
if and only if $H _{0}$ is not included in $\Phi (P)$; in particular, this
occurs in the special case of $H _{0} = P _{0}.[P, P]$, where $P _{0}$ is a
subgroup of $P$ that is not is not normal in $P$;}
\par
(iv) {\it If $P$ is not isomorphic to the quaternion group of order $8$,
then it
possesses a subgroup $P _{0}$ with the property required by (iii).}
\par
\medskip
{\it Proof.} It is well-known that $\Phi (P)$ is a normal subgroup
of $P$ including $[P, P]$, and such that $P/\Phi (P)$ is an
elementary abelian $p$-group of rank $r \ge 1$; this implies the
normality of the subgroups of $P$ including $\Phi (P)$. Recall
further that $r \ge 2$, since, otherwise, $P$ must possess exactly
one maximal subgroup, and therefore, must be nontrivial and
cyclic, in contradiction with the assumption that $P$ is
nonabelian. On the other hand, it follows from the noted
properties of $\Phi (P)$ that if $k$ is a natural number less than
$r$ and $S$ is a subset of $P$ with $k$ elements, then the
subgroup $P(S)$ of $P$ generated by the union $\Phi (P) \cup S$ is
of order dividing $\vert \Phi (P) \vert .p ^{k}$; in particular,
$P(S)$ is a proper subgroup of $P$. Since proper subgroups of $P$
are abelian, these observations show that $\Phi (P) \subseteq
Z(P)$ and $r = 2$. At the same time, the noncommutativity of $P$
ensures that $P/Z(P)$ is a noncyclic group, whence it becomes
clear that $\Phi (P) = Z(P)$. Let $h$ be an element of $P
\setminus \Phi (P)$. Then there exists an element $g \in P$, such
that the system of co-sets $\{h\Phi (P), g\Phi (P)\}$ generates
the group $P/\Phi (P)$. Using the fact that $[P, P] \subseteq \Phi
(P) = Z(P)$ and $h ^{p} \in Z(P)$, one obtains by direct
calculations that each element of $[P, P]$ is a power of the
commutator $h ^{-1}g ^{-1}hg := (h, g)$, and also, that $(h, g)
^{p} = 1$. This completes the proof of Lemma 3.2 (i). Statement
(ii) and the former part of statement (iii) of Lemma 3.2 can be
deduced from Lemma 3.2 (i). As $\Phi (P)$ consists of all
non-generators of $P$ (cf. [23, Ch. 1, Theorem 2]), the systems
$\{h, g\}$ and $\{h.[P, P], g.[P, P]\}$ generate the groups $P$
and $P/[P, P]$, respectively. Therefore, the quotient group
$P/(\Phi (P)\langle h\rangle )$ is cyclic, which implies the
latter part of Lemma 3.2 (iii). The concluding assertion of the
lemma can be obtained from the classification of Miller-Moreno
$p$-groups [31] (cf. also [37, Theorem 444]), namely, the fact
that if $P$ is not isomorphic to the quaternion group of order
$8$, then it has one of the following presentations:
\par
\medskip
$P _{1} = \langle g _{1}, h _{1}, z: g _{1} ^{p ^{m}} = h _{1} ^{p
^{n}} = z ^{p} = 1, g _{1}z = zg _{1}, h _{1}z = zh _{1}, h _{1}g
_{1}h _{1} ^{-1} = g _{1}z\rangle $, $m \ge n$ $\ge 1$ ($P _{1}$
is of order $P ^{m+n+1}$);
\par
$P _{2} = \langle g _{2}, h _{2}: g _{2} ^{p ^{m}} = h _{2} ^{p ^{n}} = 1, h _{2}g
_{2}h _{2} ^{-1} = g _{2} ^{1+p ^{m-1}}\rangle $, $ m \ge 2, n \ge 1, p
^{m+n} > 8$ ($P _{2}$ is of order $p ^{m+n}$).
\par
\medskip
Clearly, the subgroup of $P _{i}$ generated by $h _{i}$ is not
normal, for any index $i \le 2$ and any admissible pair $(m, n)$,
so Lemma 3.2 is proved.
\par
\medskip
We continue with the proof of Theorem 3.1. Suppose first that
$G(R/E)$ is not isomorphic to the quaternion group $\hbox{\Bbb Q}
_{8}$. It follows from Galois theory and Lemma 3.2 that then the
extension $R/E$ possesses an intermediate field $L$ for which the
following is true:
\par
\medskip
(3.1) (i) $R = LR _{\rm ab}$ and $L$ is not normal over $E$;
\par
(ii) The intersection $L \cap R _{\rm ab} := F$ is a cyclic extension of $E$ of
degree $[L:E]/p$ (and is not equal to $E$).
\par
\medskip
It is clear from (3.1) (ii) that $L/F$ is a cyclic extension of
degree $p$. Let $\sigma $ and $\psi $ be generators of the Galois
groups $G(L/F)$ and $G(F/E)$, respectively. Fix an element $\omega
$ of $F ^{\ast }$, denote by $\Delta $ the cyclic $F$-algebra
$(L/F, \sigma , \omega )$, and by $\bar \psi $ some embedding of
$L$ in $E (p)$ as an $E$-subalgebra, inducing $\psi $ on $F$. By
Proposition 2.5 (iii), $\psi $ is extendable to an automorphism
$\tilde \psi $ of $\Delta $ as an algebra over $E$. Observing also
that $\bar \psi (L) \neq L$ and arguing as in the proof of [9,
Lemma 3.2], one concludes that there exists an $F$-isomorphism
$\Delta \cong (L/F, \sigma , \psi (\beta )\beta ^{-1})$, for some
$\beta \in F ^{\ast }$. Hence, by [34, Sect. 15.1, Proposition b],
$\omega \beta \psi (\beta ) ^{-1}$ is an element of the norm group
$N(L/F)$. In view of [11, Lemma 3.2], this means that $\omega \in
N(R _{\rm ab}/F)$ if and only if $\omega \beta \psi (\beta ) ^{-1}
\in (N(L/F) \cap N(R _{\rm ab}/F))$. Since $F \neq E$ and $G(R/E)$
is a Miller-Moreno group, $R/F$ is an abelian extension, so it
follows from (3.1) (i) and the availability of a local $p$-class
field theory on $E$ (in the case of Br$(E) _{p} \neq \{0\}$), by
[11, Theorem 2.1], that $N(R/F) = N(L/F) \cap N(R _{\rm ab}/F)$.
Taking into consideration that $N _{E} ^{F} (\omega ) = N _{E}
^{F} (\omega \beta \psi (\beta ) ^{-1})$ one obtains from the
above results and the transitivity of norm mappings in towers of
finite extensions that $N(R/E) = N(R _{\rm ab}/E)$.
\par
Assume now that $p = 2$ and $G(R/E)$ is a quaternion group of
order $8$.  It this case, by Galois theory, $R _{\rm ab}$ is
presentable as a compositum of two different quadratic extensions
$E _{1}$ and $E _{2}$ of $E$; one also sees that $R/E _{1}$ and
$R/E _{2}$ are cyclic extensions of degree $4$.  Let $\psi _{1}$
be an $E _{1}$-automorphism of $R$ of order $4$, $\sigma _{1}$ an
$E$-automorphism of $E _{1}$ of order $2$, and $\gamma $ an
element of $R _{\rm ab} ^{\ast }$ not lying in $N(R/R _{\rm ab})$.
The field $E _{1}$ is $2$-quasilocal, which implies the existence
of a central division $E _{1}$-algebra $D$ of index $4$, such that
$D \otimes _{E _{1}} R _{\rm ab}$ is similar to the cyclic $R
_{\rm ab}$-algebra $(R/R _{\rm ab}, \psi _{1} ^{2}, \gamma )$.
Using again Proposition 2.5, one concludes that $D$ is isomorphic
to the cyclic $E _{1}$-algebra $(R/E _{1}, \psi _{1}, \rho )$, for
some $\rho \in E _{1} ^{\ast }$.  Therefore, there exists an $R
_{\rm ab}$-isomorphism $(R/R _{\rm ab}, \psi _{1} ^{2}, \gamma )
\cong (R/R _{\rm ab}, \psi _{1} ^{2}, \rho )$, and by [34, Sect.
15.1, Proposition b], $\gamma \rho ^{-1} \in N(R/R _{\rm ab})$. By
Proposition 2.5 (iii), the normality of $R$ over $E$, and the
Skolem-Noether theorem (cf. [34, Sect. 12.6]), $\sigma _{1}$ is
extendable to an automorphism $\tilde \sigma _{1}$ of $D$ as an
algebra over $E$, such that $\tilde \sigma _{1} (R) = R$.  In
addition, our assumption on $G(R/E)$ indicates that $(\psi \tilde
\sigma _{1}) (r) = (\psi _{1} ^{3}\tilde \sigma _{1}) (r)$, for
each $r \in R ^{\ast }$.  It is now easy to see that $D$ is
isomorphic to the cyclic $E _{1}$-algebras $(R/E _{1}, \psi _{1}
^{3}, \sigma _{1} (\rho ))$ and $(R/E _{1}, \psi _{1} ^{3}, \rho
^{3})$. Hence, by [34, Sect. 15.1, Proposition b], $\rho
^{3}\sigma _{1} (\rho ) ^{-1}$ lies in $N(R/E _{1})$.  Taking also
into account that $\gamma \rho ^{-1} \in N(R/R _{\rm ab})$ and $N
_{E _{1}} ^{R _{\rm ab}} (\gamma \rho ^{-1}) = \rho ^{-2}N _{E
_{1}} ^{R _{\rm ab}} (\gamma )$, one concludes that $N _{E _{1}}
^{R _{\rm ab}} (\gamma ).(\rho \sigma _{1} (\rho ) ^{-1}) \in
N(R/E _{1})$.  This result shows that $N _{E} ^{R _{\rm ab}}
(\gamma ) \in N(R/E)$, which completes the proof of Theorem 3.1 in
the special case where $R/E$ is a Galois extension.
\par
\medskip
Suppose finally that $R$ is an arbitrary extension of $E$ in $E
(p)$ of degree $p ^{m}$, for some $m \in \hbox{\Bbb N}$, and
denote by $R _{0}$ the maximal normal extension $E$ in $R$.
Proceeding by induction on $m$ and taking into account that $R
_{0} \neq E$, one obtains that now it suffices to prove Theorem
3.1 in the special case where $R \neq R _{0}$ and $R$ is abelian
over $R _{0}$, and assuming that the conclusion of the theorem is
valid for each intermediate field of $R/E$ not equal to $R$.  Our
inductive hypothesis indicates that then there exists an embedding
$\psi $ of $R$ in $E (p)$ as an $E$-subalgebra, such that $\psi
(R) \neq R$.  It is easily verified that $R _{0} = \psi (R _{0})$
and $\psi (R)/R _{0}$ is an abelian extension ; hence, $R$ and
$\psi (R)$ are abelian extensions of the intersection $R \cap \psi
(R) := R _{1}$. Observing also that $R _{1}$ is a $p$-quasilocal
field, one gets from [11, Theorem 2.1] that $R _{1} ^{\ast } =
N(R/R _{1})N(\psi (R)/R _{1})$. This, combined with the
transitivity of norm mappings in towers of finite extensions, and
with the fact that $\psi (N _{R _{0}} ^{R} (\lambda )) = N _{R
_{0}} ^{\psi (R)} (\psi (\lambda ))$, for each $\lambda \in R
^{\ast }$, implies that $N(R/E) = N(R _{1}/E)$. Since $R _{1} \neq
R, R _{0} \subseteq R _{1}$ and the normal extensions of $E$ in
$R$ are subfields of $R _{0}$, the proof of Theorem 3.1 can be
accomplished by applying the obtained result and the inductive
hypothesis.
\par
\medskip
{\bf Corollary 3.3.} {\it Let $(E, v)$ be a Henselian discrete
valued field with local class field theory, $M/E$ a finite Galois
extension with a nilpotent Galois group, and $R$ an intermediate
field of $M/E$. Then $N(R/E) = N(R _{\rm ab}/E)$.}
\par
\medskip
{\it Proof.} It suffices to consider the special case of a proper
$p$-extension $M/E$. Then $E (p) \neq E$, and by [10, Theorem
2.1], $\widehat E (p)/\widehat E$ is a $\hbox{\Bbb Z}
_{p}$-extension, where $\widehat E$ is the residue field of $(E,
v)$. If $p \neq $ char$(\widehat E)$ and $\widehat E$ does not
contain a primitive $p$-th root of unity, this means that $E
(p)/E$ is a $\hbox{\Bbb Z} _{p}$-extension with finite
subextensions inertial over $E$ (see, for example, [8, Lemma
1.1]), so our assertion becomes trivial. It follows from the
Henselian property of $v$ that if $\widehat E$ contains a
primitive $p$-th root of unity, then $E$ also contains such a
root; hence, by [11, Proposition 2.4] and the fact that $E$ admits
local $p$-class field theory, $E$ is a $p$-quasilocal field, i.e.
our statement is a special case of Theorem 3.1. Also, it has been
proved in [10, Sect. 2] that $E$ is $p$-quasilocal, provided that
char$(\widehat E) = p$ and char$(E) = 0$, so Corollary 3.3 is now
obvious.
\par
\vskip2.4truecm {\centerline{\bf 4. On nonnilpotent finite Galois
extensions}
\par
\smallskip
\centerline{\bf of strictly PQL-fields algebraic over $\hbox{\Bbb Q}$}}
\par
\medskip
Proposition 2.8 and [16, Sect. 2, Theorem 4] indicate that if $E$
is an algebraic extension of a global field $E _{0}$, $F$ is an
intermediate field of $E/E _{0}$, and for each $p \in P(E)$,
$w(p)$ is the absolute value of $F$ induced by $v(p)$, then the
groups Br$(F) _{p}$ and Br$(F _{w(p)}) _{p}$ are nontrivial.
Therefore, our research concentrates as in [13] on the study of
the following class of fields:
\par
\medskip
{\bf Definition 4.1.} Let $E _{0}$ be a global field, $\overline E
_{0}$ an algebraic closure of $E _{0}$, $F$ an extension of $E
_{0}$ in $\overline E _{0}$, $P$ a nonempty set of
prime numbers for which Br$(F) _{p} \neq \{0\}$, and $\{w(p): p
\in P\}$ a system of nontrivial absolute values of
$F$, such that Br$(F _{w(p)}) _{p} \neq \{0\}$, $p \in P$. Denote by $\Omega
(F, P, W)$ the set of intermediate fields $E$ of $\overline E/F$
with the following properties:
\par
(i) $E$ admits local class field theory and $P(E) = P$;
\par
(ii) The characteristic system $\{v(p): p \in P\}$ of $E$ can be chosen
so that $v(p)$ is a prolongation of $w(p)$, for each $p \in P$.
\par
\medskip
The main results of [13] about the set $\Omega (F, P, W)$ can be stated
as follows:
\par
\medskip
{\bf Proposition 4.2.} {\it With assumptions and notations being
as above, $\Omega (F, P, W)$ is a nonempty set, for which the
following assertions hold true:}
\par
(i) {\it Every field $E \in \Omega (P, W; F)$ possesses a unique subfield
$R(E)$ that is a minimal element of $\Omega (P, W; F)$ (with respect to
inclusion);}
\par
(ii) {\it Every minimal element $E$ of $\Omega (F, P, W)$ is an intermediate
field of a Galois extension of $F$ with a prosolvable Galois group;}
\par
(iii) {\it If $E$ is a minimal element of $\Omega (F, P, W)$, $p \in P$ and
$F _{w(p)}$ is the closure of $F$ in $E _{v(p)}$, then the degrees of the
finite extensions of $F _{w(p)}$ in $E _{v(p)}$ are not divisible by $p$.}
\par
\medskip
Proposition 4.2 plays a crucial role in the proof of the following
precise form of Theorem 1.2.
\par
\medskip
{\bf Proposition 4.3.} {\it Let $G$ be a nonnilpotent finite
group, $\overline P$ the set of all prime numbers, $w(p)$ the
normalized $p$-adic absolute value of the field $\hbox{\Bbb Q}$ of
rational numbers, $\forall p \in \overline P$, and $W = \{w(p): p
\in \overline P\}$. Then there is a field $E \in \Omega
(\hbox{\Bbb Q}, \overline P, W)$ possessing a Galois extension $M$
in $\overline {\hbox{\Bbb Q}}$, such that $G(M/E)$ is isomorphic
to $G$ and $N(M/E) \neq N(M _{\rm ab}/E)$.}
\par
\medskip
{\it Proof.} Our argument relies on several observations described
by the following four lemmas.
\par
\medskip
{\bf Lemma 4.4.} {\it Let $E$ be an algebraic strictly} PQL{\it
-extension of $\hbox{\Bbb Q}$ with $P(E) \neq \phi $, and let
$\{v(p): p \in P(E)\}$ be a characteristic system of $E$. Assume
also that $M/E$ is a finite Galois extension such that $G(M/E)$ is
nonnilpotent, each prime $p$ dividing $[M:E]$ lies in $P(E)$ and
$M _{v(p)'}/E _{v(p)}$ is a $p$-extension with $G(M _{v(p)'}/E
_{v(p)})$ isomorphic to the Sylow $p$-subgroups of $G(M/E)$, where
$v(p) ^{\prime }$ is an arbitrary absolute value of $M$ extending
$v(p)$. Then $N(M/E) \neq N(M _{\rm ab}/E)$.}
\par
\medskip
{\it Proof.} By the Burnside-Wielandt theorem, the assumption that
$G$ is nonnilpotent means that it possesses a maximal subgroup $H$
that is not normal. Let $p$ be a prime number dividing the index
$\vert G:H\vert $, $A _{p}$ the maximal $p$-extension of $E$ in
$M$, $H _{p}$ a Sylow $p$-subgroup of $H$, $G _{p}$ a Sylow
$p$-subgroup of $G(M/E)$ including $H _{p}$, and $K$, $K _{p}$ and
$M _{p}$ the intermediate fields of $M/E$ corresponding by Galois
theory to $H$, $H _{p}$ and $G _{p}$, respectively. It follows
from Galois theory and the normality of maximal subgroups of
finite $p$-groups that $A _{p} \cap K = E$, which indicates that
$(A _{p}K)/K$ is a $p$-extension with $G((A _{p}K)/K)$
isomorphic to $G(A _{p}/E)$. The extensions $(A _{p}K _{p})/K
_{p}$ and $(A _{p}M _{p})/M _{p}$ have the same property, since
the choice of $K, K _{p}$ and $M _{p}$ guarantees that the degrees
$[K _{p}:K]$ and $[M _{p}:E]$ are not divisible by $p$. This
implies that $[(A _{p}K _{p}):K _{p}] = [(A _{p}M _{p}):M _{p}] = [A
_{p}:E]$ and $[(A _{p}K _{p}):M _{p}] = [(A _{p}M _{p}):M _{p}].[K _{p}:M
_{p}]$. Thus it turns out that $A _{p}M _{p} \cap K _{p} = M _{p}$, which
means that $G _{p}$ is of greater rank as a $p$-group than $G(A _{p}/E)$, and
because of Proposition 2.10, proves Lemma 4.4.
\par
\medskip
{\bf Lemma 4.5.} {\it Let $F$ be an algebraic number field,
$\overline P$ the set of all prime numbers, $W = \{w(p): p \in
\overline P\}$ a system of nontrivial absolute values of $F$ fixed
as in Definition 4.1, $M _{0}/F$ a finite Galois extension, and
for each prime $p$ dividing $[M _{0}:F]$, let $M _{0,w(p)'}/F
_{w(p)}$ be a normal extension with a Galois group isomorphic to
the Sylow $p$-subgroups of $G(M _{0}/F)$, where $w(p) ^{\prime }$
is an arbitrary prolongation of $w(p)$ on $M _{0}$. Assume also
that $E$ is a minimal element of $\Omega (F, \overline P, W)$,
$V(E) = \{v(p): p \in \overline P\}$ is a characteristic system of
of $E$, and $M = M _{0}E$. Then $M/E$ is a Galois extension
satisfying the conditions of Lemma 4.4, and the Galois groups
$G(M/E)$ and $G(M _{0}/E _{0})$ are canonically isomorphic.}
\par
\medskip
{\it Proof.} Denote by $P(M _{0}/F)$ the set of prime divisors of
$[M _{0}:F]$. The minimality of $E$ and Proposition 4.2 (iii)
imply that $M _{0,w(p)'} \otimes _{F _{w(p)}} E _{v(p)}$ is a field
isomorphic to $M _{v(p)'}$ over $E _{v(p)}$, where $v(p) ^{\prime
}$ is a prolongation of $v(p)$ on $M _{0}$, for each $p \in P(M
_{0}/E _{0})$. This, combined with the fact that the groups
$G(M/E)$ and $G(M _{v(p)'}/E _{v(p)})$ embed in $G(M _{0}/F)$ and
$G(M/E)$, respectively, and with the condition on the extension $M
_{0,w(p)'}/F _{w(p)}$, proves that $G(M/E)$ is isomorphic to $G(M
_{0}/F)$.
\par
\medskip
{\bf Lemma 4.6.} {\it Let $M/E$ be a Galois extension with a
Galois group $G$ embeddable in the symmetric group $S _{n}$, for
some $n \in \hbox{\Bbb N}$. Then there exists a polynomial $f(X)
\in E[X]$ of degree $n$ with a root field (over $E$) equal to
$M$.}
\par
\medskip
{\it Proof.} Denote by $s$ the number of $G$-orbits of the set
$\{1,..., n\}$, fix a system $\{g _{j}: j = 1,..., s\}$ of
representatives of these orbits, and for each index $j$, let $U
_{j}$ be the intermediate field of $M/E$ corresponding by Galois
theory to the stabilizer Stab$_{G} (g _{j}) := G _{j}$. It is
easily verified that $[U _{j}:E] = \vert G:G _{j}\vert $, $\sum
_{j=1} ^{s} \vert G:G _{j}\vert = n$, and $\cap _{j=1} ^{s} V _{j}
= \{1\}$, where $V _{j}$ is the intersection of the subgroups of
$G$ conjugate to $G _{j}$, $\forall j \in \{1,..., s\}$.
Therefore, one can take as $f(X)$ the product $\prod _{j=1} ^{s} f
_{j}(X)$, where $f _{j} (X)$ is the minimal polynomial over $E$ of
any primitive element of $U _{j}$ over $E$, for each $j$.
\par
\medskip
{\bf Lemma 4.7.} {\it Let $E _{0}$ be an algebraic number field,
$n$ an integer number greater than one, and $P _{n}$ the set of
prime numbers $\le n$. Assume that $E _{0}$ possesses, for each $p
\in P _{n}$, an absolute value $w(p)$ for which the completion $E
_{0,w(p)}$ admits a Galois extension $\widetilde M _{p}$ with
$G(\widetilde M _{p}/E _{0,w(p)})$ isomorphic to the Sylow
$p$-subgroups of the symmetric group $S _{n}$. Then there exists a
Galois extension $M _{0}$ of $E _{0}$ with $G(M _{0}/E _{0})$
isomorphic to $S _{n}$, and such that the completion $M
_{0,w(p)'}$ is $E _{0,w(p)}$-isomorphic to $\widetilde M _{p}$,
for each $p \in P _{n}$, and any prolongation $w(p) ^{\prime }$ of
$w(p)$ on $M$.}
\par
\medskip
{\it Proof.} It follows from Lemma 4.6 and the assumptions of the
present lemma that $\widetilde M _{p}$ is a root field over $E
_{0,w(p)}$ of a separable polynomial $f _{p} (X) = X ^{n} + \sum
_{j=1} ^{n} c _{p,j}X ^{n-j} \in $ $E _{0,w(p)} [X]$, for each $p
\in P _{n}$. Since the absolute values $w(p): p \in P _{n}$ are
pairwise nonequivalent in the case of $n > 2$, the weak
approximation theorem (cf. [26, Ch. XII, Sect. 1]) and the density
of $E _{0}$ in $E _{0,w(p)}$ ensure, for each real positive number
$\varepsilon $, the existence of a polynomial $g _{\varepsilon }
(X) \in E _{0} [X]$ equal to $X ^{n} + \sum _{j=0} ^{n-1} b
_{\varepsilon ,j}X ^{n-j}$, and such that $w(p) (b _{\varepsilon
,j} - c _{p,j}) < \varepsilon $, for every $p \in P _{n}$. This
enables one to deduce from Krasner$^{,}$s lemma (cf. [27, Ch. II,
Proposition 3]) that if $\varepsilon $ is sufficiently small, then
the quotient rings $E _{0,w(p)} [X]/f _{p} (X)E _{0,w(p)} [X]$ and
$E _{0,w(p)} [X]/g _{\varepsilon } (X)E _{0,w(p)} [X]$ are
isomorphic as $E _{0,w(p)}$-algebras, which implies that
$\widetilde M _{p}$ is a root field of $g _{\varepsilon } (X)$
over $E _{0,w(p)}$, $\forall p \in P _{n}$. When this occurs, it
becomes clear from Galois theory and the obtained result that the
root field of $g _{\varepsilon } (X)$ over $E _{0}$ is a normal
extension of $E _{0}$ with a Galois group $G _{\varepsilon }$ of
order divisible by $n!$. As $G _{\varepsilon }$ obviously embeds
in $S _{n}$, this means that $G _{\varepsilon } \cong S _{n}$, so
Lemma 4.7 is proved.
\par
\medskip
We are now in a position to prove Proposition 4.3. Retaining
assumptions and notations in accordance with Lemma 4.5, note that
every intermediate field $K$ of $M _{0}/E _{0}$ possesses a system
$\{\nu (p): p \in P(M _{0}/E _{0})\}$ of absolute values, such
that $\nu (p)$ is a prolongation of $w(p)$ and $M _{\nu (p)'}/K
_{\nu (p)}$ is a normal extension with a Galois group isomorphic
to the Sylow $p$-subgroups of $G(M _{0}/K)$, $\forall p \in P(M
_{0}/K)$. To show this, take a prime $p \in P(M _{0}/K)$, fix a
Sylow $p$-subgroup $P$ of $G(M _{0}/K)$ as well as a Sylow
$p$-subgroup $P _{0}$ of $G(M _{0}/E _{0})$ including $P$, and
denote by $F _{0}$ and $F$ the extensions of $E _{0}$ in $M _{0}$
corresponding by Galois theory to $P _{0}$ and $P$, respectively.
The local behaviour of $M _{0}/E _{0}$ at $w(p) ^{\prime }/w(p)$
implies the existence of a prolongation $\omega _{0} (p)$ of
$w(p)$ on $F _{0}$, such that $F _{0,\omega _{0} (p)}$ is a
completion of $E _{0}$ with respect to $w(p)$; moreover, it
becomes clear that $\omega _{0} (p)$ is uniquely extendable to an
absolute value $\omega (p)$ of $M _{0}$ (cf. [6, Ch. II, Theorem
10.2]). Observing that $[F:F _{0}] = \vert P _{0}:P\vert $ and $p$
does not divide $[F _{0}:K]$, one concludes that the absolute
value $\nu (p)$ of $K$ induced by $\omega (p)$ has the required
property. Since every finite group of order $n$ is embeddable in
the symmetric group $S _{n}$, for each $n \in \hbox{\Bbb N}$
(Cayley's theorem), the obtained result and the previous three
lemmas indicate that Proposition 4.3 will be proved, if we show
the existence of an algebraic number field $E _{0}$ satisfying the
conditions of Lemma 4.7.
\par
Fix a natural number $n > 1$ as well as an odd integer $m > n!$,
suppose that $P _{n}$ is defined as in Lemma 4.7, put $\tilde n =
\prod _{p \in P _{n}} p$, and denote by $\Phi _{0}$ the extension
of $\hbox{\Bbb Q}$ in $\overline {\hbox{\Bbb Q}}$ obtained by
adjoining a root of the polynomial $X ^{m} - \tilde n$. Also, let
$\Gamma _{s} = \hbox{\Bbb Q} (\delta _{s} + \delta _{s} ^{-1})$,
$\Phi _{s} = \Phi _{0} (\delta _{s} + \delta _{s} ^{-1})$, where
$\delta _{s}$ is a primitive $2 ^{s}$-th root of unity in
$\overline {\hbox{\Bbb Q}}$, $\forall s \in {\Bbb N}$, and $\Phi
_{\infty } = \cup _{s=1} ^{\infty } \Phi _{s}$. The choice of
$\Phi _{0}$ indicates that the $p$-adic absolute value of
$\hbox{\Bbb Q}$ is uniquely extendable to an absolute value $w
_{0} (p)$ of $\Phi _{0}$, for each $p \in P _{n}$; one obtains
similarly that the $2$-adic absolute value of $\hbox{\Bbb Q}$ has
a unique prolongation $w _{s} (2)$ on $\Phi _{s}$, $\forall s \in
\hbox{\Bbb N}$ (cf. [6, Ch. I, Theorem 6.1]), and also, a unique
prolongation $w _{\infty }$ on $\Phi _{\infty }$. Furthermore, our
argument proves that $\Phi _{s,w _{s}(2)}/\hbox{\Bbb Q} _{2}: s
\in {\Bbb N}$ and $\Phi _{0,w _{0}(p)}/{\Bbb Q} _{p}: p \in P
_{n}, p > 2$, are totally ramified extensions of degrees $2
^{s}.m$ and $m$, respectively. We first show that $\Phi _{0,w
_{0}(p)}$ admits a Galois extension with a Galois group isomorphic
to the Sylow $p$-subgroups of $S _{n}$, for each $p \in (P _{n}
\setminus \{2\})$. Note that $\Phi _{0,w _{0}(p)}$ does not
contain a primitive $p$-th root of unity. This follows from the
fact that $m$ is odd whereas $p - 1$ equals the degree of the
extension of $\hbox{\Bbb Q} _{p}$ obtained by adjoining a
primitive $p$-th root of unity (cf. [18, Ch. IV, (1.3)]). Hence,
by the Shafarevich theorem [44] (cf. also [42, Ch. II, Theorem
3]), the Galois group of the maximal $p$-extension of $\Phi _{0,w
_{0}(p)}$ is a free pro-$p$-group of rank $m + 1$. In view of
Galois theory, this means that a finite $p$-group is realizable as
a Galois group of a $p$-extension of $\Phi _{0,w _{0}(p)}$ if and
only if it is of rank at most equal to $m + 1$. The obtained
result, combined with the fact that $m > n!$ and the ranks of the
$p$-subgroups of $S _{n}$ are less than $n!$, proves our
assertion. Taking now into consideration that $\Phi _{s}/\Phi
_{0}$ is a cyclic extension of degree $2 ^{s}$, one obtains by
applying [6, Ch. II, Theorem 10.2] and [26, Ch. IX, Proposition
11] that $\Phi _{s,w _{s}(p)}$ is a cyclic extension of $\Phi
_{0,w _{0}(p)}$ of degree dividing $2 ^{s}$, for each absolute
value $w _{s} (p)$ of $\Phi _{s}$ extending $w _{0} (p)$. It is
therefore clear from Galois theory that Proposition 4.3 will be
proved, if we show that $\Phi _{s,w _{s}(2)}$ admits a normal
extension with a Galois group isomorphic to the Sylow
$2$-subgroups of $S _{n}$, for every sufficient large index $s$.
Identifying $\Phi _{0,w _{0}(2)}$ with the closure of $\Phi _{0}$
in $\Phi _{\infty ,w _{\infty }}$, one obtains from the uniqueness
of the prolongation $w _{\infty }/w _{0} (2)$ that $\Phi _{\infty
} \cap \Phi _{0,w _{0}(2)} = \Phi _{0}$ and the compositum $\Phi
_{\infty }\Phi _{0,w _{0}(2)} := \Phi _{\infty ,2}$ is a
$\hbox{\Bbb Z} _{2}$-extension of $\Phi _{0,w _{0}(2)}$. This
implies that Br$(\Phi _{\infty ,2}) _{2} = \{0\}$ and $\Phi
_{\infty ,2} (2) \neq \Phi _{\infty ,2}$, which means that $G(\Phi
_{\infty ,2} (2)/\Phi _{\infty ,2})$ is a free pro-$2$-group of
countably infinite rank (cf. [42, Ch. II, 5.6, Theorem 4 and Lemma
3] and [48, p. 725]). Hence, finite $2$-groups are realizable as
Galois groups of normal extensions of $\Phi _{\infty ,2}$. In
particular, there exists a $2$-extension $T _{\infty ,2}$ of $\Phi
_{\infty ,2}$ with $G(T _{\infty ,2}/\Phi _{\infty ,2})$
isomorphic to the Sylow $2$-subgroups of $S _{n}$, so it follows
from [12, (1.3)] that one can find an index $\tilde s$ and a
Galois extension $T _{\tilde s,2}$ of $\Phi _{\tilde s,w _{\tilde
s}(2)}$ in $T _{\infty ,2}$, such that $T _{\tilde s,2} \otimes
_{\Phi _{\tilde s,w _{\tilde s}(2)}} \Phi _{\infty ,2}$ is
isomorphic to $T _{\infty ,2}$ as an algebra over $\Phi _{\tilde
s,w _{\tilde s}(2)}$. In view of the general properties of tensor
products (cf. [34, Sect. 9.4, Corollary a]), this implies that if
$s$ is an integer $\ge \tilde s$, then the $\Phi _{s,w
_{s}(2)}$-algebra $T _{\tilde s,2} \otimes _{\Phi _{\tilde s,w
_{\tilde s}(2)}} \Phi _{s,w _{s}(2)} := T _{s,2}$ is a field, and
more precisely, a Galois extension of $\Phi _{s,w _{s}(2)}$ with
$G(T _{s,2}/\Phi _{s,w _{s}(2)})$ isomorphic to the Sylow
$2$-subgroups of $S _{n}$. Furthermore, in this case, the field
$\Phi _{s} := E _{0}$ and its absolute values $w _{s} (p)$, $p \in
P _{n}$, satisfy the conditions of Lemma 4.7, which completes the
proof of Proposition 4.3.
\par
\vskip0.6truecm \centerline{\bf 5. Proof of Theorem 1.1 (ii)}
\par
\medskip
Let $E$ be a field, $R/E$ a finite separable extension, and for
each prime $p$, let $R _{{\rm ab},p}$ be the maximal abelian
$p$-extension of $E$ in $R$, $\rho _{p}$ the greatest integer
dividing $[R:E]$ and not divisible by $p$, and $N _{p} (R/E)$ the
set of those elements $u _{p} \in E ^{\ast }$, for which the
co-set $u _{p}N(R/E)$ is a $p$-element of the  group $E ^{\ast
}/N(R/E)$. Clearly, $u ^{\rho _{p}} \in N _{p} (R/E)$, for every
$u \in E ^{\ast }$. Observing also that $u ^{\rho _{p}} \in N(R
_{{\rm ab},p}/E)$ whenever $u \in N(R _{\rm ab}/E)$ and $p$ is
prime, one obtains that Theorem 1.1 (ii) can be deduced from the
following result.
\par
\medskip
{\bf Theorem 5.1.} {\it Assume that $E$ is a quasilocal field,
such that the natural homomorphism of} Br$(E)$ {\it into} Br$(L)$
{\it maps} Br$(E) _{p}$ {\it surjectively on} Br$(L) _{p}${\it ,
for some prime number $p$ and every finite extension $L$ of $E$.
Then $N(R/E)$ includes as a subgroup the intersection $N(R _{{\rm
ab},p}/E) \cap N _{p} (R/E)$, for each finite extension $R$ of $E$
in $E _{\rm sep}$.}
\par
\medskip
The rest of this Section is devoted to the proof of Theorem 5.1.
If Br$(E) _{p} = \{0\}$, our assumptions ensure that Br$(L) _{p} =
\{0\}$, for every finite extension $L$ of $E$, which reduces our
assertion to a special case of Lemma 2.4. Assuming further that
Br$(E) _{p} \neq \{0\}$, and $\hbox{\Bbb F} _{p}$ is a field with
$p$ elements (identifying it with the prime subfield of $E$, in
the case of char$(E) = p$), we first prove Theorem 5.1 in the
special case where the Galois group $G(\widetilde R/E)$ of the
normal closure $\widetilde R$ of $R$ in $E _{\rm sep}$ over $E$ is
solvable. The main part of our argument is presented by the
following three lemmas.
\par
\medskip
{\bf Lemma 5.2.} {\it Let $E$ be a field and $p$ a prime number
satisfying the conditions of Theorem 5.1, and let $M/E$ be a
Galois extension with $G(M/E)$ possessing the following two
properties:}
\par
(i) $G(M/E)$ {\it is nonabelian and isomorphic to a semidirect
product $E _{p;k} \times C _{\pi }$ of an elementary abelian
$p$-group of order $p ^{k}$ by a group $C _{\pi }$ of prime order
$\pi $ not equal to $p$, and $k$ is the minimal positive integer
for which $p ^{k}$ is congruent to $1$ modulo $\pi $;}
\par
(ii) $E _{p;k}$ {\it is a minimal normal subgroup of $G(M/E)$.
\par
Then the norm group $N(M/E _{1})$ includes $E ^{\ast }$, where $E
_{1}$ is the intermediate field of $M/E$ corresponding by Galois
theory to $E _{p;k}$.}
\par
\medskip
{\it Proof.} Our assumptions indicate that $E _{1}/E$ is a cyclic
extension of degree $\pi $, and under the additional hypothesis
that Br$(E) _{p} \neq \{0\}$, this means that Br$(E _{1}) _{p}
\neq \{0\}$ (see [34, Sect. 13.4]). Therefore, by [11, Theorem
2.1], $E _{1}$ admits local $p$-class field theory, so it is
sufficient to show that $E ^{\ast } \subseteq N(M _{1}/E _{1})$,
for every cyclic extension $M _{1}$ of $E _{1}$ in $M$. Suppose
first that $E$ contains a primitive $p$-th root of unity or
char$(E) = p$, and fix an $E$-automorphism $\psi $ of $E _{1}$ of
order $\pi $. As $G(M/E _{1})$ is an elementary abelian $p$-group
of rank $k$, Kummer$^{,}$s theory and the Artin-Schreier theorem
imply the existence of a subset $S = \{\rho _{j}: j = 1,..., k\}$
of $E _{1}$, such that the root field over $E _{1}$ of the
polynomial set $\{f _{j} (X) = X ^{p} - uX - \rho _{j}: j = 1,...,
k\}$ equals $M$, where $u = 1$, if char$(E) = p$, and $u = 0$,
otherwise. For each index $j$, denote by $z _{j}$ the element
$\psi (u _{j})u _{j} ^{-1}$ in case $E$ contains a primitive
$p$-th root of unity, and put $z _{j} = \psi (u _{j}) - u _{j}$,
if char$(E) = p$. Note that $M$ is a root field over $E _{1}$ of
the set of polynomials $\{g _{j} (X) = X ^{p} - uX - z _{j}: j =
1,..., k\}$. This can be deduced from the following two
statements:
\par
\medskip
(5.1) (i) If char$(E) = p$, $r(E _{1}) = \{\lambda ^{p} - \lambda
: \lambda \in E _{1}\}$, and $M(E _{1})$ is the additive subgroup
of $E$ generated by the union $S \cup r(E _{1})$, then $r(E _{1})$
and $M(E _{1})$ are $\psi $-invariant, regarded as vector spaces
over $\hbox{\Bbb F} _{p}$; moreover, the linear operator of the
quotient space $M(E _{1})/r(E _{1})$, induced by $\psi - id _{E
_{1}}$ is an isomorphism;
\par
(ii) If $E$ contains a primitive $p$-th root of unity, $M(E _{1})$
is the multiplicative subgroup of $E _{1} ^{\ast }$ generated by
the union $S \cup E _{1} ^{\ast p}$, and the mapping $\psi _{1}: E
_{1} ^{\ast }/E _{1} ^{\ast p} \to E _{1} ^{\ast }/E _{1} ^{\ast
p}$ is defined by the rule $\psi _{1} (\alpha  E _{1} ^{\ast p}) =
\psi (\alpha )\alpha ^{-1}E _{1} ^{\ast p}: \alpha \in E _{1}
^{\ast }$, then $\psi _{1}$ is a linear operator of $E _{1} ^{\ast
}/E _{1} ^{\ast p}$ (regarded as a vector space over $\hbox{\Bbb
F} _{p}$), $M(E _{1})/E _{1} ^{\ast p}$ is a $k$-dimensional $\psi
_{1}$-invariant subspace of $E _{1} ^{\ast }/E _{1} ^{\ast p}$,
and the linear operator of $M(E _{1})/E _{1} ^{\ast p}$ induced by
$\psi _{1}$ is an isomorphism.
\par
\medskip
Most of the assertions of (5.1) are well-known. One should,
possibly, only note here that the concluding parts of (5.1) (i)
and (5.1) (ii) follow from the fact that $G(M/E _{1})$ is the
unique normal proper subgroup of $G(M/E)$, and by Galois theory,
this means that $E _{1}$ is the unique normal proper extension of
$E$ in $M$. The obtained result implies the nonexistence of a
cyclic extension of $E$ in $M$ of degree $p$, which enables one to
deduce from Kummer's theory and the Artin-Schreier theorem the
triviality of the kernels of the considered linear operators. Thus
our argument leads to the conclusion that the discussed special
case of Lemma 5.2 will be proved, if we establish the validity of
the following two statements, for each index $j$:
\par
\medskip
(5.2) (i) If $E$ contains a primitive $p$-th root of unity
$\varepsilon $, and $c$ is an element of $E ^{\ast }$, then the
symbol $E _{1}$-algebra $A _{\varepsilon } (z _{j}, c; E _{1})$ is
trivial;
\par
(ii) If char$(E) = p$ and $c \in E ^{\ast }$, then the $p$-symbol
$E$-algebra $E[z _{j}, c)$ is trivial.
\par
\medskip
Denote by $D _{j}$ the symbol $p$-algebra $E _{1} [\rho _{j}, c)$,
if char$(E) = p$, and the symbol $E _{1}$-algebra $A _{\varepsilon
} (\rho _{j}, c; E _{1})$, in case $E _{1}$ contains a primitive
$p$-th root of unity $\varepsilon $. It follows from the
assumptions of Theorem 5.1 that $D _{j}$ is isomorphic over $E
_{1}$ to $\Delta _{j} \otimes _{E} E _{1}$, for some central
division $E$-algebra $\Delta _{j}$. In view of the Skolem-Noether
theorem, this implies the extendability of $\psi $ to an
automorphism $\bar \psi $ of $D _{j}$, regarded as an algebra over
$E$. Thus it becomes clear that $D _{j}$ is $E _{1}$-isomorphic to
$E _{1} [\psi (\rho _{j}), c)$ or $A _{\varepsilon } (\psi (\rho
_{j}), c; E _{1})$ depending on whether or not char$(E) = p$.
Applying now the general properties of local symbols (cf. [43, Ch.
XIV, Propositions 4 and 11]), one proves (5.2).
\par
It remains for us to prove Lemma 5.2, assuming that $p \neq $
char$(E)$ and $E$ does not contain a primitive $p$-th root of
unity. Let $\varepsilon $ be such a root in $M _{\rm sep}$. It is
easily verified that if $E(\varepsilon ) \cap E _{1} = E$, then
$M(\varepsilon )/E(\varepsilon )$ is a Galois extension, such that
$G((M(\varepsilon )/E(\varepsilon ))$ is canonically isomorphic to
$G(M/E)$. Since $E(\varepsilon )$ and $p$ satisfy the conditions
of the lemma, our considerations prove in this case that
$E(\varepsilon ) ^{\ast } \subseteq N(M(\varepsilon )/E _{1}
(\varepsilon ))$. Hence, by Lemma 2.1, applied to the triple $(E
_{1}, M, E _{1} (\varepsilon ))$ instead of $(E, L _{1}, L _{2})$,
we have $E ^{\ast } \subseteq N(M/E _{1})$, which reduces the
proof of Lemma 5.2 to the special case in which $E _{1}$ is an
intermediate field of $E(\varepsilon )/E$. Fix a generator
$\varphi $ of $G(E(\varepsilon )/E)$, and an integer number $s$ so
that $\varphi (\varepsilon ) = \varepsilon ^{s}$. Observing that
$M/E$ is a noncyclic Galois extension of degree $p\pi $, one
obtains from (2.2) and the cyclicity of $M$ over $E _{1}$ that
$M(\varepsilon )$ is generated over $E(\varepsilon )$ by a $p$-th
root of an element $\rho $ of $E(\varepsilon )$ with the property
that $\varphi (\rho )\rho ^{-s'} \in E(\varepsilon ) ^{\ast p}$,
where $s ^{\prime }$ is a positive integer such that $s ^{\prime
\pi } \equiv s ^{\pi }$(mod $p$) and $s ^{\prime } \not\equiv
s$(mod $p$). It is therefore clear from (2.3) and [34, Sect. 15.1,
Proposition b] that $A _{\varepsilon } (\rho , c; E(\varepsilon
))$ is isomorphic to the matrix $E(\varepsilon )$-algebra $M _{p}
(E(\varepsilon ))$, $\forall c \in E ^{\ast }$. One also sees that
$E ^{\ast } \subseteq N(M(\varepsilon )/E(\varepsilon ))$. As
$[M:E _{1}] = p$ and $[E(\varepsilon ):E _{1}]$ divides $(p -
1)/\pi $, Lemma 2.1 ensures now that $E ^{\ast } \subseteq N(M/E
_{1})$, so Lemma 5.2 is proved.
\par
\medskip
{\bf Lemma 5.3.} {\it Assuming that $E$ is a quasilocal field
whose finite extensions satisfy the conditions of Theorem 5.1, for
a given prime number $p$, suppose that $M/E$ is a finite Galois
extension, such that $G(M/E)$ is a solvable group. Then $N _{p}
(M/E) \cap N(M _{{\rm ab},p}/E)$ is a subgroup of $N(M/E)$.}
\par
\medskip
{\it Proof.} It is clearly sufficient to prove the lemma under the
hypothesis that $N(M ^{\prime }/E ^{\prime })$ includes $N _{p} (M
^{\prime }/E ^{\prime }) \cap N(M ^{\prime } _{{\rm ab},p}/E
^{\prime })$, provided that $E ^{\prime }$ and $p$ satisfy the
conditions of Theorem 5.1, and $M ^{\prime }/E ^{\prime }$ is a
Galois extension with a solvable Galois group of order less than
$[M:E]$. As in the proof of Theorem 1.1 (i), we first show that
then one may assume further that $G(M/E)$ is a Miller-Moreno
group. Our argument relies on the fact that the class of fields
satisfying the conditions of Theorem 5.1 is closed under the
formation of finite extensions. Note that if $G(M/E)$ is not a
Miller-Moreno group, then it possesses a nonabelian subgroup $H$
whose commutator subgroup $[H, H]$ is normal in $G(M/E)$. Indeed,
one can take as $H$ the commutator subgroup $[G(M/E), G(M/E)]$, in
case $G(M/E)$ is not metabelian, and suppose that $H$ is any
nonabelian maximal subgroup of $G(M/E)$, otherwise. Denote by $F$
and $L$ the intermediate fields of $M/E$ corresponding to $H$ and
$[H, H]$, respectively. Our choice of $H$ and Galois theory
indicate that $L$ is a Galois extension of $E$ including $M _{\rm
ab}$, and such that $E \neq L \neq M$, so our additional
hypothesis and Lemma 2.2 lead to the conclusion that $N _{p} (L/E)
\cap N(M _{{\rm ab},p}/E) = N _{p} (L/E) \cap N(M _{\rm ab}/E)
\subseteq N(L/E)$ and $N _{p} (M/F) \cap N(L/F) \subseteq N(M/F)$.
Let now $\mu $ be an element of $N _{p} (M/E) \cap N(M _{\rm
ab}/E)$, and $\lambda \in L ^{\ast }$ a solution to the norm
equation $N _{E} ^{L} (X) = \mu $. Then one can find an integer
$k$ not divisible by $p$, and such that $N _{F} ^{L} (\lambda )
^{k} \in N _{p} (M/F)$. It is therefore clear that $N _{F} ^{L}
(\lambda ) ^{k} \in N(M/F)$ and $\mu ^{k} \in N(M/E)$. As $\mu \in
N _{p} (M/E)$, this implies that $\mu \in N(M/E)$, which yields
the desired reduction. In view of Theorem 1.1 (i) and the
elementary properties of norm mappings, one may also assume that
$G(M/E)$ is a nonnilpotent Miller-Moreno group, such that $p$
divides the index $\vert G:[G, G]\vert $. By the classification of
these groups [31] (cf. also [37, Theorem 445]), this means that
$G(M/E)$ has the following structure:
\par
\medskip
(5.3) (i) $G(M/E)$ is isomorphic to a semi-direct product $E
_{p;k} \times C _{\pi ^{n}}$ of $E _{p;k}$ by a cyclic group $C
_{\pi ^{n}}$ of order $\pi ^{n}$, for some different prime numbers
$p$ and $\pi $, where $k$ satisfies condition (i) of Lemma 5.2;
\par
(ii) $E _{p;k}$ is a minimal normal subgroup of $G(M/E)$, and the
centre of $G(M/E)$ equals the subgroup $C _{\pi ^{n-1}}$ of $C
_{\pi ^{n}}$ of order $\pi ^{n-1}$.
\par
\medskip
It follows from (5.2) and Galois theory that the extension, say,
$E _{n}$ of $E$ in $M$ corresponding to $E _{p;k}$ is cyclic of
degree $\pi ^{n}$. This indicates that $N _{E _{n}} ^{M} (\eta
_{n}) = \eta _{n} ^{p ^{k}}$, for every $\eta _{n} \in E _{n}$,
and thereby, proves that $c ^{p ^{k}} \in N(M/E)$, in case $c \in
N(E _{n}/E)$. We show in this case that $c ^{\pi ^{n}} \in
N(M/E)$. By Lemma 5.2, if $n = 1$, then $E _{1} ^{\ast }$ contains
an element $\xi $ of norm $c$ over $E _{1}$, which means that $N
_{E} ^{M} (\xi ) = c ^{\pi }$. Suppose now that $n \ge 2$, put
$\tilde \pi = \pi ^{n-1}$, denote by $C _{\tilde \pi }$ the
subgroup of $G(M/E)$ of order $\tilde \pi $, and let $M ^{\prime
}$ and $E ^{\prime }$ be the intermediate fields of $M/E$
corresponding by Galois theory to the subgroups $C _{\tilde \pi }$
and $E _{p;k}C _{\tilde \pi }$ of $G(M/E)$, respectively. It is
easily seen that $M ^{\prime }/E$ is a Galois extension with $G(M
^{\prime }/E)$ satisfying the conditions of Lemma 5.2, and $E
^{\prime }/E$ is a cyclic extension of degree $\pi $. This ensures
that $c ^{\pi } \in N(M ^{\prime }/E)$. Also, it becomes clear
that $M = M ^{\prime }E _{n}$, $M ^{\prime } \cap E _{n} = E
^{\prime }$, and $N _{M'} ^{M} (m ^{\prime }) = m ^{\prime \tilde
\pi }$: $m ^{\prime } \in M ^{\prime }$. These observations show
that $c ^{\pi ^{n}} \in N(M/E)$. Since $c ^{p ^{k}} \in N(M/E)$,
g.c.d.$(p, \pi ) = 1$, and $E _{n} \subseteq M _{\rm ab} \subset
M$, the obtained result implies the inclusions $N(M _{\rm ab}/E)
\subseteq N(E _{n}/E) \subseteq N(M/E) \subseteq N(M _{\rm
ab}/E)$. Evidently, these inclusions are equalities, so Lemma 5.3
is proved.
\par
\medskip
{\bf Lemma 5.4.} {\it Retaining assumptions and notations as in
Theorem 5.1, suppose that $M/E$ is a Galois extension with a
solvable Galois group $G(M/E)$, and $R$ is an intermediate field
of $M/E$, such that $[R:E]$ is a power of $p$. Then $N(R/E) = N(R
_{\rm ab}/E)$.}
\par
\medskip
{\it Proof.} Arguing by induction on $[M:E]$, one obtains from the
conditions of Theorem 5.1 that it is sufficient to prove the
lemma, assuming in addition that $N(R _{1}/E _{1}) = N(R ^{\prime
}/E _{1})$, whenever $E _{1}$ and $R _{1}$ are intermediate fields
of $M/E$, such that $E _{1} \neq E$, $E _{1} \subseteq R _{1}$,
$[R _{1}:E _{1}]$ is a power of $p$, and $R ^{\prime }$ is the
maximal abelian extension of $E _{1}$ in $R _{1}$. Suppose first
that $R _{\rm ab} \neq E$. Then the inductive hypothesis, applied to
the the pair $(E _{1}, R _{1}) = (R _{\rm ab}, R)$, gives $N(R/E) =
N(R ^{\prime }/E)$, and since $R ^{\prime }$ is a subfield of the
maximal $p$-extension $M _{p}$ of $E$ in $M$, this enables one to
obtain from Theorem 1.1 (i) that $N(R ^{\prime }/E) = N(R
_{\rm ab}/E)$.
\par
It remains to be seen that $N(R/E) = E ^{\ast }$ in the special
case of $R _{\rm ab} = E$. Our argument relies on the fact that $E
^{\ast }/N(M _{\rm ab}/E)$ is a group of exponent dividing $[M
_{\rm ab}:E]$. Therefore, if $M _{p} = E$, then this exponent is
not divisible by $p$. In view of the inclusion $N(M/E) \subseteq
N(R/E)$, $E ^{\ast }/N(R/E)$ is canonically isomorphic to a
quotient group of $E ^{\ast }/N(M/E)$, so the condition $M _{p} =
E$ ensures that the exponent $e(R/E)$ of $E ^{\ast }/N(R/E)$ is
also relatively prime to $p$. As $e(R/E)$ divides $[R:E]$, this
proves that $N(R/E) = E ^{\ast }$.
\par
Assume now that $R _{\rm ab} = E$ and $M _{p} \neq E$, denote by
$F _{1}$ the maximal abelian extension of $E$ in $M _{p}$, and by
$F _{2}$ the intermediate field of $M/E$ corresponding by Galois
theory to some Sylow $p$-subgroup of $G(M/E)$. Put $R _{1} = RF
_{1}$, $R _{2} = RF _{2}$ and $F _{3} = F _{1}F _{2}$. It follows
from Galois theory and the equality $R _{\rm ab} = E$ that the
compositum $RM _{p}$ is a Galois extension of $R$ with $G((RM
_{p})/R)$ canonically isomorphic to $G(M _{p}/E)$; in addition, it
becomes clear that $R _{1}$ is the maximal abelian extension of
$R$ in $RM _{p}$. Thus it turns out that $[R _{1}:R] = [F
_{1}:E]$, which means that $[R _{1}:E] = [R:E].[F _{1}:E]$.
Observing that $[F _{2}:E]$ is not divisible by $p$, one also sees
that $[R _{2}:F _{2}] = [R:E]$, $[(R _{1}F _{2}):F _{2}] = [R
_{1}:E]$ and $[(RF _{3}):F _{2}] = [R _{2}:F _{2}].[F _{3}:F
_{2}]$. The concluding equality implies that $R _{2} \cap F _{3} =
F _{2}$. As $F _{2}$ admits local $p$-class field theory, in case
Br$(E) _{p} \neq \{0\}$, this leads to the conclusion that $N(R
_{2}/E)N(F _{3}/E) = N(F _{2}/E)$. Note also that $N(F _{1}/E) =
N(R _{1}/E)$. Indeed, it follows from Galois theory and the
definition of $M _{p}$ that $M _{p}$ does not admit proper
$p$-extensions in $M$, and by the inductive hypothesis, this
yields $N((RM _{p})/M _{p}) = M _{p} ^{\ast }$. Hence, by Theorem
1.1 (i) and the general properties of norm mappings, we have
$N((RM _{p})/E) = N(M _{p}/E) = N(F _{1}/E)$. At the same time,
since $R _{1}$ is the maximal abelian extension of $R$ in $RM
_{p}$, it turns out that $N((RM _{p})/R) = N(R _{1}/R)$, which
implies that $N((RM _{p})/E) = N(R _{1}/E) = N(F _{1}/E)$, as
claimed. The obtained results and the inclusions $N(R _{2}/E)
\subseteq N(R/E)$ and $N(F _{3}/E) \subseteq N(F _{1}/E)$,
indicate that $N(F _{2}/E)$ is a subgroup of $N(R/E)N(F _{1}/E)$
$= N(R/E)N(R _{1}/E) = N(R/E)$. As $E ^{\ast }/N(R/E)$ and $E
^{\ast }/N(F _{2}/E)$ are groups of finite relatively prime
exponents, this means that $N(R/E) = E ^{\ast }$, so the proof of
Lemma 5.4 is complete.
\par
\medskip
We are now in a position to prove Theorem 5.1 in the special case
of a solvable group $G(M/E)$. It is clearly sufficient to
establish our assertion under the additional hypothesis that $N
_{p} (R _{1}/E _{1})$ and $N(R _{1}/E _{1})$ are related in
accordance with Theorem 5.1, whenever $E _{1}$ and $R _{1}$ are
extensions of $E$ in $R$ and $M$, respectively, such that $E _{1}
\neq E$ and $E _{1} \subseteq R _{1}$. Suppose that $R \neq E$,
put $\Phi = R _{{\rm ab},p}$, if $R _{\rm ab} \neq E$, and denote
by $\Phi $ some extension of $E$ in $R$ of primary degree,
otherwise (the existence of $\Phi $ in the latter case follows
from Galois theory and the well-known fact that maximal subgroups
of solvable finite groups are of primary indices). Also, let
$\alpha $ be an element of $N _{p} (M/E) \cap N(M _{{\rm
ab},p}/E)$, $\Phi ^{\prime }$  the maximal abelian $p$-extension
of $\Phi $ in $R$, and $M ^{\prime }$ the compositum $\Phi M
_{{\rm ab},p}$. It is not difficult to see that $\Phi ^{\prime }
\cap M ^{\prime } = \Phi $. Using the fact that $\Phi $ is a field
with local $p$-class field theory, if Br$(\Phi ) _{p} \neq \{0\}$,
and applying Lemma 2.4, if Br$(\Phi ) _{p} = \{0\}$, one obtains
further that $\Phi ^{\ast } = N(\Phi ^{\prime }/\Phi )N(M ^{\prime
}/\Phi )$. At the same time, Lemma 5.4 and the choice of $\Phi $
ensure the existence of an element $\xi \in \Phi $ of norm $\alpha
$ over $E$. Observing also that there is a natural number $k$ not
divisible by $p$, for which $\Phi ^{\ast k} \subseteq N _{p}
(M/\Phi )$, one obtains from the inductive hypothesis, the
equality $\Phi ^{\ast } = N(\Phi ^{\prime }/\Phi )N(M ^{\prime
}/\Phi )$, and the inclusions $N(M/\Phi ) \subseteq N(R/\Phi )$,
$N(M ^{\prime }/E) \subseteq N(M _{\rm ab},p/E)$, that $\xi ^{k}
\in N(R/\Phi )N(M ^{\prime }/\Phi )$ and $\alpha ^{k} \in N(R/E)(N
_{p} (M/E) \cap N(M _{{\rm ab},p}/E))$. Hence, by Lemma 5.3,
$\alpha ^{k} \in N(R/E)$, and since $\alpha \in N _{p} (R/E)$ and
$p$ does not divide $k$, this means that $\alpha \in N(R/E)$,
which proves Theorem 5.1 in the special case where $G(M/E)$ is
solvable. In order to establish the theorem in full generality, we
need the following lemma.
\par
\medskip
{\bf Lemma 5.5.} {\it Assuming that $E$ and $p$ satisfy the
conditions of Theorem 5.1, suppose that $R$ is an intermediate
field of a finite Galois extension $M/E$, such that $G(M/E)$
equals the commutator subgroup $[G(M/E), G(M/E)]$. Then $N _{p}
(R/E) \subseteq N(R/E)$.}
\par
\medskip
{\it Proof.} It is clearly sufficient to consider only the special
case of $R = M \neq E$ (and Br$(E) _{p} \neq \{0\}$). Denote by $E
_{p}$ be the intermediate field of $M/E$ corresponding by Galois
theory to some Sylow $p$-subgroup of $G(M/E)$. Then $p$ does not
divide the degree $[E _{p}:E] := m _{p}$, so the condition Br$(E)
_{p} \neq \{0\}$ guarantees that Br$(E _{p}) _{p} \neq \{0\}$. We
first show that $E ^{\ast } \subseteq N(M/E _{p})$, assuming
additionally that char$(E) = p$ or $E$ contains a primitive root
of unity of degree $[M:E _{p}]$. As $E$ is a quasilocal field, the
nontriviality of Br$(E _{p}) _{p}$ ensures that $E _{p}$ admits
local $p$-class field theory. Hence, by Theorem 1.1 (i), it is
sufficient to prove the inclusion $E ^{\ast } \subseteq N(L/E
_{p})$, for an arbitrary cyclic extension $L$ of $E _{p}$ in $M$.
By [34, Sect. 15.1, Proposition b], this is equivalent to the
assertion that the cyclic $E _{p}$-algebra $(L/E _{p}, \sigma ,
c)$ is isomorphic to the matrix $E _{p}$-algebra $M _{n} (E
_{p})$, where $c \in E ^{\ast }$, $n = [L:E _{p}]$ and $\sigma $
is an $E _{p}$-automorphism of $L$ of order $n$. Since g.c.d.$([E
_{p}:E], p) = 1$, the surjectivity of the natural homomorphism of
Br$(E) _{p}$ into Br$(E _{p}) _{p}$ implies that the corestriction
homomorphism cor$_{E _{p}/E}$: Br$(E _{p}) \to $Br$(E)$ induces an
isomorphism of Br$(E _{p}) _{p}$ on Br$(E) _{p}$ (cf. [46, Theorem
2.5]). Applying the projection formula (cf. [28, Proposition 3
(i)] and [46, Theorem 3.2]), as well as Kummer$^{,}$s theory and
its analogue for finite abelian $p$-extensions over a field of
prime characteristic $p$, due to Witt (see, for example, [24, Ch.
7, Sect. 3]), one obtains that cor$_{E _{p}/E}$ maps the
similarity class $[(L/E _{p}, \sigma , c)]$ into $[(\widetilde
L/E, \tilde \sigma , c)]$, for some cyclic $p$-extension
$\widetilde L$ of $E$ in $M$. As $G(M/E) = [G(M/E), G(M/E)]$,
these observations show that $\widetilde L = E$ and $[(\widetilde
L/E, \tilde \sigma , c)] = 0$ in Br$(E)$. Furthermore, it becomes
clear that $[(L/E _{p}, \sigma , c)] = 0$ in Br$(E _{p})$, i.e. $c
\in N(L/E _{p})$, which proves the inclusion $E ^{\ast } \subseteq
N(M/E _{p})$. Since $N _{E} ^{E _{p}} (c) = c ^{m _{p}}$, one also
sees that $c ^{m _{p}} \in N(M/E)$: $c \in E ^{\ast }$, in the
special case where $p = $ char$(E)$ or $E$ contains a primitive
root of unity of degree $[M:E _{p}]$.
\par
Suppose now that $p \neq $ char$(E)$, fix a primitive root of
unity $\varepsilon \in M _{\rm sep}$ of degree $[M:E _{p}]$, and
put $\Phi (\varepsilon ) = \Phi ^{\prime }$, for every
intermediate field $\Phi $ of $M/E$, and $H ^{m _{p}} = \{h ^{m
_{p}}: h \in H\}$, for each subgroup $H$ of $M ^{\prime \ast }$.
As $E ^{\prime }/E$ is an abelian extension, our assumption on
$G(M/E)$ ensures that $E ^{\prime } \cap M = E$, and by Galois
theory, this means that $M ^{\prime }/E ^{\prime }$ is a Galois
extension with $G(M ^{\prime }/E ^{\prime })$ canonically
isomorphic to $G(M/E)$. Thus it becomes clear from the preceeding
considerations that $E ^{\prime \ast m _{p}} \subseteq N(M
^{\prime }/E ^{\prime })$ and $N(E ^{\prime }/E) ^{m _{p}}
\subseteq N(M ^{\prime }/E) \subseteq N(M/E)$. Our argument also
shows that $M \cap E _{p} ^{\prime } = E _{p}$, and since $E _{p}$
is $p$-quasilocal, it enables one to deduce from [11, Theorem
2.1], Theorem 1.1 (i) and Lemma 2.2 that $N(M/E _{p})N(E _{p}
^{\prime }/E _{p}) = E _{p} ^{\ast }$. Hence, by the transitivity
of norm mappings in towers of finite extensions, $N(M/E)N(E _{p}
^{\prime }/E) = N(E _{p}/E)$. These observations prove the
inclusions $E ^{\ast m _{p} ^{2}} \subseteq N(E _{p}/E) ^{m _{p}}
\subseteq N(M/E) ^{m _{p}}.N(E ^{\prime }/E) ^{m _{p}} \subseteq
N(M/E)$. This, combined with the fact that $p$ does not divide $m
_{p}$ and $[M:E]$ divides the exponent of the group $E ^{\ast
}/N(M/E)$, indicates that $E ^{\ast m _{p}} \subseteq N(M/E)$, and
so completes the proof of Lemma 5.5.
\par
\medskip
It is now easy to complete the proof of Theorem 5.1. Assume that
$M _{0}$ is the maximal Galois extension of $E$ in $M$ with a
solvable Galois group, and also, that $\mu _{p}$, $m _{p}$ and
$\rho _{p}$ are the maximal integer numbers not divisible by $p$
and dividing $[M _{0}:E]$, $[M:E]$ and $[R:E]$, respectively.
Applying Lemma 5.3 to $M _{0}/E$ and Lemma 5.5 to $M/M _{0}$, one
obtains that $E ^{\ast \mu _{p}} \subseteq N(M _{0}/E)$ and $M
_{0} ^{\ast \bar m _{p}} \subseteq N(M/M _{0})$, where $\bar m
_{p} = m _{p}/\mu _{p}$. Hence, by the norm transitivity identity
$N _{E} ^{M} = N _{E} ^{M _{0}} \circ N _{M _{0}} ^{M}$, we have
$E ^{\ast m _{p}} \subseteq N(M _{0}/E) ^{\bar m _{p}} \subseteq
N(M/E)$. Since $E ^{\ast [R:E]} \subseteq N(R/E)$, $N(M/E)
\subseteq N(R/E)$ and g.c.d.$(m _{p}, [R:E]) = \rho _{p}$, this
means that $E ^{\ast \rho _{p}} \subseteq N(R/E)$, so Theorem 5.1
is proved.
\par
\medskip
{\bf Remark 5.6.} (i) The fulfillment of the conditions of Theorem
1.1 (ii) is guaranteed, if $E$ is a field with local class field
theory in the sense of Neukirch-Perlis [33], i.e. if the triple
$(G _{E}, \{G(E _{\rm sep}/F), F \in \Sigma \}, E _{\rm sep}
^{\ast })$ is an Artin-Tate class formation (cf. [3, Ch. XIV]),
where $\Sigma $ is the set of finite extensions of $E$ in $E _{\rm
sep}$. When this occurs, the assertion of Theorem 1.1 (ii) is
contained in [3, Ch. XIV, Theorem 7], and examples of this kind
are given by $p$-adically closed or Henselian discrete valued
field with quasifinite residue field (see [36, Theorem 3.1 and
Lemma 2.9] and [43, Ch. XIII, Proposition 6]). Let us note without
going into details that the class of fields satisfying the
conditions of Theorem 1.1 (ii) includes properly the one studied
in [33].
\par
(ii) The question of whether the absolute Galois groups of
quasilocal fields are prosolvable seems to be open. Its answer is
affirmative in the special cases considered in Sections 6 and 7
(cf. [7, Proposition 3.1] and [8, Sect. 3]). It is worth noting in
this connection that the proof of Theorems 5.1 and 1.1 (ii) in the
special case where $G(M/E)$ is a solvable group bears an explicit
field-theoretic character.
\par
\vskip0.6truecm
{\centerline{\bf 6. Henselian discrete valued fields whose finite}
\par
\smallskip
\noindent
\centerline{\bf extensions are strictly primarily quasilocal}}
\par
\medskip
In this Section, we use Theorem 5.1 for describing the norm groups
of finite separable extensions of the fields pointed out in its
title, and for characterizing those of them, whose absolute Galois
groups and finite extensions have the properties required by
Theorem 1.3. Throughout the Section, $\overline P$ is the set of
prime numbers, and for each field $E$, $P _{0} (E)$ is the subset
of those $p \in \overline P$, for which $E$ contains a primitive
$p$-th root of unity, or else, $p = $ char$(E)$. Also, we denote
by $P _{1} (E)$ the subset of $\overline P \setminus P _{0} (E)$,
of those numbers $p ^{\prime }$, for which $E ^{\ast } \neq E
^{\ast p'}$, and put $P _{2} (E) = \overline P \setminus (P _{0}
(E) \cup P _{1} (E))$. Every finite extension $L$ of a field $K$
with a Henselian valuation $v$ is considered with its valuation
extending $v$, this prolongation is also denoted by $v$, $U(L)$
and $e(L/K)$ denote the multiplicative group of the valuation ring
of $(L, v)$, and the ramification index of $L/K$, respectively,
and $U(L) ^{\nu } = \{\lambda ^{\nu }: \lambda \in U(L)\}$, for
each $\nu \in {\Bbb N}$. Our starting point is the following
statement (proved in [10]):
\par
\medskip
(6.1) With assumptions being as above, if $v$ is discrete, then the
following conditions are equivalent:
\par
(i) Finite extensions of $K$ are strictly PQL-fields;
\par
(ii) The residue field $\widehat K$ of $(K, v)$ is perfect, the
absolute Galois group $G _{\widehat K}$ is metabelian of
cohomological $p$-dimension cd$_{p}(G _{K}) = 1$, for each prime
$p$, and $P _{0} (\widetilde L) \subseteq P(\widetilde L)$, for
every finite extension $\widetilde L$ of $\widehat K$.
\par
When these conditions are in force, $P _{j} (K) \setminus
\{$char$(\widehat K)\} = P _{j} (\widehat K) \setminus
\{$char$(\widehat K)\} $: $j = 0, 1$, the quotient group
$\widetilde L ^{\ast }/\widetilde L ^{\ast p ^{\nu }}$ is cyclic
of order $p ^{\nu }$, for every $\nu \in \hbox{\Bbb N}$ and each
$p \in (P _{0} (\widehat K) \cup P _{1} (\widehat K))$, $p \neq $
char$(\widehat K)$. Also, Br$(L) _{p}$ is isomorphic to the
quasicyclic $p$-group $\hbox{\Bbb Z} (p ^{\infty })$, for every
finite extension $L$ of $K$, and each $p \in P(\widehat L)$.
\par
\medskip
The main result of this Section, stated below, sheds light on the
norm groups of finite separable extensions of $K$, provided that
it satisfies the equivalent conditions in (6.1). Applied to the
special case where $P(\widehat K) = \overline P$, it yields the
norm limitation theorem for Henselian discrete valued fields with
quasifinite residue fields.
\par
\medskip
{\bf Theorem 6.1.} {\it Assume that $(K, v)$ is a Henselian
discrete valued field whose finite extensions are strictly}
PQL{\it , and let $R$ be a finite extension of $K$ in $K _{\rm
sep}$. Then $R/K$ possesses an intermediate field $R _{1}$, for
which the following is true:}
\par
(i) {\it The sets of  prime divisors of $e(R _{1}/K)$, $[\widehat
R _{1}:\widehat K]$ and $[\widehat R:\widehat R _{1}]$ are
included in $P _{1} (\widehat K)$, $\overline P \setminus
P(\widehat K)$ and $P(\widehat K)$, respectively;}
\par
(ii) $N(R/K) = N((R _{\rm ab}R _{1})/K)$ {\it and $K ^{\ast
}/N(R/K)$ is isomorphic to the direct sum $G(R _{\rm ab}/K) \times
(K ^{\ast }/N(R _{1}/K))$; in particular, $K ^{\ast }/N(R/K)$ is
of order $[R _{\rm ab}:K][R _{1}:K]$.}
\par
\medskip
{\it Proof.} Let $R ^{\prime }$ be the maximal inertial extension
of $K$ in $R$, i.e. the inertial lift of $\widehat R$ in $R$ over
$K$ (cf. [22, Theorems 2.8 and 2.9]). Note first that $R ^{\prime
}$ contains as a subfield an extension of $K$ of degree $n _{0}$,
for each $n _{0} \in \hbox{\Bbb N}$ dividing $[R ^{\prime }:K]$.
Indeed, the validity of (6.1) (ii) indicates that $G _{\widehat
K}$ is metabelian and cd$_{p} (G _{\widehat K}) = 1$: $p \in
\overline P$, and by [7, Lemma 1.2], this means that the Sylow
pro-$p$-subgroups of $G _{\widehat K}$ are isomorphic to
$\hbox{\Bbb Z} _{p}$, $\forall p \in \overline P$. It is therefore
clear that the Sylow subgroups of the Galois groups of finite
Galois extensions of $\widehat K$ are cyclic. Let now $M$ be the
normal closure of $R ^{\prime }$ in $K _{\rm sep}$ over $K$. It is
well-known that then $M$ is inertial over $K$ and the Galois
groups $G(M/K)$ and $G(\widehat M/\widehat K)$ are isomorphic (cf.
[22, page 135]). This enables one to deduce our assertion from
Galois theory and the following lemma.
\par
\medskip
{\bf Lemma 6.2.} {\it Assume that $G$ is a nontrivial finite group
whose Sylow subgroups are cyclic, $H$ is a subgroup of $G$ of
order $n$, and $n _{1}$ is a positive integer dividing the order
$o(G)$ of $G$ and divisible by $n$. Then $G$ possesses a subgroup
$H _{1}$ of order $n _{1}$, such that $H \subseteq H _{1}$.}
\par
\medskip
{\it Proof.} Our assumptions show that $G$ is a supersolvable
group, and therefore, it has a normal Sylow $p$-subgroup $G _{p}$,
as well as a subgroup $A _{p}$ isomorphic to $G/G _{p}$, where $p$
is the greatest prime divisor of the order $o(G)$ of $G$ (cf. [23,
Ch. 7, Sect. 1, Theorem 4, and Sect. 2]). Proceeding by induction
on $o(G)$, one obtains from this result (and the supersolvability
of subgroups of $G$) that $G$ possesses a subgroup $\widetilde H
_{1}$ of order $n _{1}$. Now the conclusion of the lemma follows
from the fact [38] (see also [45, Theorem 18.7]) that $H$ is
conjugate in $G$ to a subgroup of $\widetilde H _{1}$.
\par
\medskip
Let now $R _{1} ^{\prime }$ be the maximal tamely ramified
extension of $K$ in $R$, $[R _{1} ^{\prime }:R ^{\prime }] = n$,
$P$ the set of prime numbers dividing $[R:K]$, and for each $p \in
P$, let $f(p)$ and $g(p)$ be the greatest nonnegative integers for
which $p ^{f(p)} \vert [R:K]$ and $p ^{g(p)} \vert [R ^{\prime
}:K]$. As noted above, Lemma 6.2 indicates that there is an
extension $R _{p}$ of $K$ in $R ^{\prime }$ of degree $p ^{g(p)}$,
$\forall p \in P$. Observing that $\alpha \in U(R ^{\prime })
^{n}$, provided that $\alpha \in R ^{\prime }$ and $v(\alpha - 1)
> 0$, one obtains from [27, Ch. II, Proposition 12] that $R _{1}
^{\prime } = R ^{\prime } (\theta )$, where $\theta $ is an $n$-th
root of $\pi \rho $, for a suitably chosen element $\rho \in U(R
^{\prime })$. Suppose now that $p \in (P _{0} (\widehat K) \cup P
_{1} (\widehat K))$ and $p \neq $ char$(\widehat K)$. Since $p$
does not divide $[R ^{\prime }:R _{p}]$, then the concluding
assertion of (6.1) implies the existence of an element $\rho _{p}
\in U(R _{p})$, such that $\rho _{p}\rho ^{-1}$ is a $p
^{(f(p)-g(p))}$-th power in $U(R ^{\prime })$. Therefore, the
binomial $X ^{p ^{(f(p)-g(p))}} - \pi \rho _{p}$ has a root
$\theta _{p} \in R _{1} ^{\prime }$. Summing up these results, one
proves the following:
\par
\medskip
(6.2) For each $p \in P \cap (P _{0} (\widehat K) \cup P _{1}
(\widehat K))$, $p \neq $ char$(\widehat K)$, there exists an
extension $T _{p}$ of $K$ in $R _{1} ^{\prime }$ of degree $p
^{f(p)}$; moreover, if $p \in P _{0} (\widehat K)$, then the
normal closure of $T _{p}$ in $K _{\rm sep}$ over $K$ is a
$p$-extension.
\par
\medskip
Denote by $R _{1}$ the compositum of the fields $R _{p}$: $p \in
(P _{2} (\widehat K) \setminus P(\widehat K))$, and $T _{p}$: $p
\in P _{1} (\widehat K)$, and put $L = R ^{\prime }R _{1}$, $T = R
_{\rm ab}R _{1}$, $L ^{\prime } = R ^{\prime } (\theta ^{\mu })$
and $T ^{\prime } = T(\theta ^{\mu })$, where $\mu $ is the
greatest integer dividing $[R _{1} ^{\prime }:R _{1}]$ and not
divisible by any element of $\overline P \setminus P(\widehat K)$.
It is easily verified that $R _{1}$ has the properties required by
Theorem 6.1 (i), $R ^{\prime } \subseteq T$, $R _{1} \subseteq L
^{\prime } \subseteq R _{1} ^{\prime }$, and $T ^{\prime } = R
_{\rm ab}L ^{\prime }$. As g.c.d.$([R _{\rm ab}:K], [R _{1}:K]) =
1$, and by (6.1), Br$(\widehat K) _{p} \cong \hbox{\Bbb Z} (p
^{\infty })$, $p \in P(\widehat K)$, Lemma 2.1 and [11, Theorem 2]
indicate that $K ^{\ast }/N(T/K) \cong G(R _{\rm ab}/K) \times (K
^{\ast }/N(R _{1}/K))$. One also sees that $e(T/K) =$ $e(R _{\rm
ab}/K)e(R _{1}/K)$, and $L ^{\prime }/L$ and $T ^{\prime }/T$ are
tamely totally ramified extensions of degree $\prod p
^{f(p)-g(p)}$, where $p$ ranges over the elements of $P \cap P
_{2} (\widehat K)$. It remains to be proved that $N(T/K) = N(R/K)$
and $K ^{\ast }/N(R _{1}/K)$ is a (cyclic or a $2$-generated)
group of order $[R _{1}:K]$. Our argument is based on the
following two statements:
\par
\medskip
(6.3) (i) The natural homomorphism of Br$(K)$ into Br$(Y)$ maps
Br$(K) _{p}$ surjectively on Br$(Y) _{p}$, for every finite
extension $Y$ of $K$ in $K _{\rm sep}$, and each $p \in P(\widehat
Y)$;
\par
(ii) $N(T ^{\prime }/K) = N(T/K)$ and $r\pi ^{[R':K]} \in N(R/K)$,
for some element $r \in U(K)$.
\par
\medskip
Statement (6.3) (i) is implied by the final assertion of (6.1) and
the well-known fact (cf. [34, Sects. 13.4 and 14.4]) that the
relative Brauer group Br$(Y/E)$ is of exponent dividing $[Y:E]$.
The rest of the proof of (6.3) relies on the fact that $R ^{\prime
}$ is the maximal inertial extension of $K$ in $R$. In particular,
$R$ is totally ramified over $R ^{\prime }$, which means that $U(R
^{\prime })$ contains an element $\rho $, such that $\rho \pi \in
N(R/R ^{\prime })$. Therefore, the latter part of (6.3) (ii)
applies to the element $r = N _{K} ^{R'} (\rho )$. In view of
(6.1) and Galois cohomology (cf. [42, Ch. II, Proposition 6 (b)]),
we have $N(\widehat R/\widehat K) = \widehat K ^{\ast }$, so it
follows from the Henselian property of $v$ that $N(R ^{\prime }/K)
= U(K)\langle r\pi ^{[R':K]}\rangle $. Since $R ^{\prime }
\subseteq T \subseteq T ^{\prime } \subseteq R _{1} ^{\prime }$,
$T$ and $T ^{\prime }$ are tamely and totally ramified over $R
^{\prime }$, and these observations show that $N(T/K) = U(K)
^{e(T/K)}\langle r\pi ^{[R':K]}\rangle $ and $N(T ^{\prime }/K) =
U(K) ^{e(T'/K)}\langle r\pi ^{[R':K]}\rangle $. As proved above,
$[T ^{\prime }:T]$ is not divisible by any $p \in (P _{0}
(\widehat K) \cup P _{1} (\widehat K))$, whereas $e(T/K) =$ $e(R
_{\rm ab}/K)e(R _{1}/K)$, so it turns out that g.c.d.$([T ^{\prime
}:T], e(T/K)) = 1$, $U(K) ^{e(T/K)}$ $= U(K) ^{e(T'/K)}$ and
$N(T/K) = N(T ^{\prime }/K)$. Arguing in a similar manner, one
obtains that $N(R _{1}/K) = U(K) ^{e(R _{1}/K)}\langle r\pi ^{[R
_{0}:K]}\rangle $, where $R _{0} = R ^{\prime } \cap R _{1}$.
Since prime divisors of $[R _{1}:K]$ lie in $P _{1} (\widehat K)$
(and $e(R _{1}/K)$ divides $[R _{1}:K]$), the concluding assertion
of (6.1) implies that $U(K)/U(K) ^{e(R _{1}/K)}$ is a cyclic group
of order $e(R _{1}/K)$. Thus the required properties of $K ^{\ast
}/N(R _{1}/K)$ become obvious.
\par
We turn to the proof of the equality $N(R/K) = N(T/K)$. The
inclusion $N(R/K)$ $\subseteq N(T/K)$ is obvious, so we prove the
inverse one. Consider an arbitrary element $\beta $ of $U(K) \cap
N(T/K)$, put $[R:T ^{\prime }] = m$, and for each $p \in (P \cap
P(\widehat K))$, denote by $R _{{\rm ab},p}$ the maximal abelian
$p$-extension of $K$ in $R$, and by $m _{p}$ the greatest integer
dividing $[R:K]$ and not divisible by $p$. It follows from the
inclusion $N(T/K) \subseteq N(R _{{\rm ab},p}/K)$, statement (6.3)
(i) and Theorem 5.1 that $\beta ^{m _{p}} \in N(R/K)$, $\forall p
\in (P \cap P (\widehat K))$. At the same time, the equality
$N(T/K) = N(T ^{\prime }/K)$ implies that $\beta ^{m} \in N(R/K)$.
Observing now that prime divisors of $m$ lie in $P(\widehat K)$,
one obtains that g.c.d.$\{m, m _{p}: p \in (P \cap P(\widehat
K))\} = 1$, and therefore, $\beta \in N(R/K)$. Since $N(T/K) =
U(K) ^{e(T/K)}\langle r\pi ^{[R':K]}\rangle $ and $N(R/K) = U(K)
^{e(R/K)}\langle r\pi ^{[R':K]}\rangle $, this means that $N(T/K)
\subseteq N(R/K)$, so Theorem 6.1 is proved.
\par
\medskip
{\bf Corollary 6.3.} {\it Assume that $K$ satisfies the conditions
of Theorem 6.1, and $R$ is a finite extension of $K$ in $K _{\rm
sep}$, such that $N(R/K) = N(\Phi (R)/K)$, for some abelian finite
extension $\Phi (R)$ of $K$ in $K _{\rm sep}$. Then $\Phi (R) = R
_{\rm ab}$.}
\par
\medskip
{\it Proof.} This follows at once from Theorem 6.1.
\par
\medskip
{\bf Remark 6.4.} Let $(K, v)$ be a Henselian discrete valued
field satisfying the conditions of Theorem 6.1.
\par
(i) It is easily deduced from (6.1) that if $n$ is a positive
integer not divisible by char$(\widehat K)$, then $K ^{\ast n}$ is
an open subgroup of $K ^{\ast }$ (with respect to the topology
induced by $v$), such that $K ^{\ast }/K ^{\ast n}$ is isomorphic
to a direct product of a cyclic group of order $n$ by a cyclic
group of order $n _{0}n _{1}$, where $n _{0}$ is the greatest
divisor of $n$, for which $\widehat K$ contains a primitive $n
_{0}$-th root of unity, and $n _{1}$ is the greatest divisor of
$n$ not divisible by any $t \in (\overline P \setminus P _{1})$.
In addition, it is not difficult to see that every subgroup $U$ of
$K ^{\ast }$ of index $n$ includes $K ^{\ast n}$ and equals
$N(R(U)/K)$, for some finite extension $R(U)$ of $K$ of degree
$n$. Furthermore, it becomes clear that if $e$ is the exponent of
$K ^{\ast }/U$, then $K ^{\ast e} \subseteq U$,  $n/e$ divides
g.c.d.$\{e, n _{0}\}$ and $K ^{\ast }/U$ is presentable as a
direct sum of a cyclic group of order $e$ by such a group of order
$n/e$. When char$(\widehat K) = 0$, this fully describes norm
groups of $K$.
\par
(ii) Suppose now that char$(\widehat K) = p > 0$ and $L$ is a
finite separable extension of $K$. Then it follows from Theorem
6.1 and Lemmas 2.2 and 2.2 that $N(L/K) = N(L _{{\rm ab},p}$ $\cap
N(L _{0}/K)$, for some intermediate field $L _{0}$ of $L/K$ of
degree $[L _{0}:K]$ not divisible by $p$. Hence, by Hazewinkel's
existence theorem [20] (see also [17, 3.5 and 3.7]) concerning
totally ramified finite abelian $p$-extensions of $K$, $N(L/K)$ is
an open subgroup of $K ^{\ast }$ of finite index. As a matter of
fact, Hazewinkel's theorem and the equality $N(L/K) = N(L _{{\rm
ab},p}/K) \cap N(L _{0}/K)$ allow one to obtain a satisfactory
inner characterization of norm groups of $K$ (see also [41] for
the special case in which char$(K) = p$ and $K$ is complete). At
the same time, it should be noted that if $\widehat K$ is of
infinite cardinality $\kappa $, then the set Nr$(K)$ of these
groups and the one of finite abelian $p$-extensions of $K$ in $K
_{\rm sep}$ are of cardinality $\kappa $ whereas the set Op$(K)$
of open subgroups of $K ^{\ast }$ of finite indices is of
cardinality $2 ^{\kappa }$. This is obtained in the spirit of the
proof of [49, Part IV, Propositions 3 and 4] or of the Corollary
of [18, Ch. V, (3.6)] (see also the proof of [10, Lemma 2.3 (i)]
in the case of char$(K) = 0$). Thus it becomes clear that Nr$(K)
\neq $ Op$(K)$ unless $\widehat K$ is a finite field.
\par
\medskip
{\bf Proposition 6.5.} {\it For a Henselian discrete valued field
$(K, v)$, the following condititions are equivalent:}
\par
(i) $G _{K}$ {\it and the finite extensions of $K$ have the
properties required by Theorem 1.3;}
\par
(ii) char$(\widehat K) = 0${\it , $P _{0} (K) = P(K) \neq
\overline P$ and $P _{1} (K) = \overline P \setminus P _{0} (K)$.
\par
When this occurs, every finite extension $R$ of $K$ in $K _{\rm
sep}$ is presentable as a compositum $R = R _{0}R _{1}$, where $R
_{1}$ has the properties described in Theorem 6.1 (i) and (ii),
and $R _{0}$ is an intermediate field of $R/K$ of degree $[R
_{0}:K] = [R:R _{1}]$. Moreover, the Galois group of the normal
closure $\widetilde R$ of $R$ in $K _{\rm sep}$ over $K$ is
nilpotent if and only if $R = R _{0}$.}
\par
\medskip
{\it Proof.} The implication (ii)$\to $(i) and the concluding
assertions of Proposition 6.5 are obtained in a straightforward
manner from Theorem 6.1 and statement (6.2), so we assume further
that condition (i) is in force. As in the proof of Theorem 6.1,
let $\pi $ be a generator of the maximal ideal of the valuation
ring of $(K, v)$. It is easily deduced from Theorem 6.1 that if $p
\in (\overline P \setminus (P _{0} (K) \cup P _{1} (K))$, or $p =
$ char$(\widehat K)$ and $p \not\in P _{0} (K)$, then the root
field, say $M _{\pi }$, of the binomial $X ^{p} - \pi $ satisfies
the equality $N(M _{\pi }/K) = N(M _{\pi ,{\rm ab}}/K)$. At the
same time, it follows from (2.2) that $G(M _{\pi }/K)$ is
nonabelian and isomorphic to a semidirect product of a group of
order $p$ by a cyclic group of order dividing $p - 1$. This
indicates that $G(M _{\pi }/K)$ is nonnilpotent. The obtained
results contradict condition (i), and thereby, prove that $P _{0}
(K) \cup P _{1} (K)$ includes $P \setminus \{$char$(\widehat
K))\}$. Our argument, together with the concluding part of (6.1),
also shows that $P _{0} (\widehat K) \cup P _{1} (\widehat K) =
\overline P$, which implies that $\widehat K$ is infinite. It
remains to be seen that char$(\widehat K) = 0$. Suppose that
char$(\widehat K) = q > 0$ and $q \in P _{0} (K)$. Then condition
(i) and statement (6.1) imply the existence of a primitive $p$-th
root of unity in $\widehat K$, for at least one prime number $p
\neq q$. In addition, since $\widehat K$ is infinite, it becomes
clear that there exists a cyclic inertial extension $L _{p}$ of
$K$ in $K _{\rm sep}$ of degree $p$. One also obtains from Galois
theory (cf. [26, Ch. VIII, Theorem 20]) and the Henselian property
of $v$ that $L _{p}$ possesses a normal basis $B _{p}$ over $K$,
such that $B _{p} \subset U(L _{p})$. Denote by $B _{p} ^{\prime
}$ the polynomial set $\{X ^{q} - X - b\pi ^{-1}: b \in B _{p}\}$,
if char$(K) = q$, and put $B _{p} ^{\prime } = \{X ^{q} - (1 +
b\pi ): b \in B _{p}\}$, in the mixed-characteristic case. It
follows from the Artin-Schreier theorem, Capelli's criterion (cf.
[26, Ch. VIII, Sect. 9]) and the Henselian property of the
prolongation of $v$ on $L _{p}$ that $B _{p} ^{\prime }$ consists
of irreducible polynomials over $L _{p}$. Furthermore, one obtains
from Kummer's theory (and the assumption that $q \in P _{0} (K)$)
that the root field $L _{p} ^{\prime }$ of $B _{p} ^{\prime }$
over $L _{p}$ is a Galois extension of $K$ of degree $q ^{p}p$. It
follows from the definition of $L _{p} ^{\prime }$ that the Sylow
$q$-subgroup $G(L _{p} ^{\prime }/L _{p})$ of $G(L _{p} ^{\prime
}/K)$, is normal and elementary abelian. At the same time, it is
clear from the choice of $B _{p}$ that $G(L _{p} ^{\prime }/L
_{p})$ possesses maximal subgroups that are not normal in $G(L
_{p} ^{\prime }/E)$. These properties of $G(L _{p} ^{\prime }/L
_{p})$ indicate that $G(L _{p} ^{\prime }/E)$ is nonnilpotent. On
the other hand, since $q \in P(\widehat K)$, Theorem 6.1 shows
that $N(L _{p} ^{\prime }/K) = N(L _{p,{\rm ab}} ^{\prime }/K)$.
Thus the hypothesis that char$(\widehat K) \neq 0$ leads to a
contradiction with condition (i), which completes the proof of
Proposition 6.5.
\par
\medskip
{\bf Corollary 6.6.} {\it Let $(K, v)$ be a Henselian discrete
valued field satisfying the conditions of Proposition 6.5,
$\varepsilon _{p}$ a primitive $p$-th root of unity in $K _{\rm
sep}$, for each $p \in \overline P$, and $[K(\varepsilon _{p}):K]
= \gamma _{p}$, in case $p \in (\overline P \setminus P(\widehat
K))$. Assume also that $L$ is a finite extension of $K$ in $K
_{\rm sep}$. Then some of the following assertions holds true:}
\par
(i) {\it $G _{L}$ and finite extensions of $L$ have the properties
required by Theorem 1.3;}
\par
(ii) $G _{L}$ {\it is pronilpotent.
\par
The latter occurs if and only if the set $\Gamma (K) = \{\gamma _{p}$: $p \in
(\overline P \setminus P(\widehat K))\}$ is bounded and $L$ contains as a
subfield the inertial extension of $K$ in $\overline K$ of degree equal to the
least common multiple of the elements of $\Gamma (K)$.}
\par
\medskip
{\it Proof.} The fulfillment of the conditions of Proposition 6.5
guarantees that $P _{0} (\widehat L)$ $\cup P _{1} (\widehat L) =
\overline P$. Applying Galois theory and (6.2), one also obtains
that $P _{0} (\widehat L) = \overline P$ if and only if $G _{L}$
is pronilpotent, which completes the proof of the corollary.
\par
\medskip
Our next result supplements Theorem 6.1, and combined with
Proposition 6.5, proves Theorem 1.3.
\par
\medskip
{\bf Proposition 6.7.} {\it Let $P _{0}$, $P _{1}$, $P _{2}$ and
$P$ be subsets of the set $\overline P$ of prime numbers, such
that $P _{0} \cup P _{1} \cup P _{2} = \overline P$, $2 \in P
_{0}$, $P _{i} \cap P _{j} = \phi $: $0 \le i < j \le 2$, and $P
_{0} \subseteq P \subseteq (P _{0} \cup P _{2})$. For each $p \in
(P _{1} \cup P _{2})$, let $\gamma _{p}$ be an integer $\ge 2$
dividing $p - 1$ and not divisible by any element of $\overline P
\setminus P$. Assume also that $\gamma _{p} \ge 3$, in case $p \in
(P _{2} \setminus P)$. Then there exists a Henselian discrete
valued field $(K, v)$ satisfying the following conditions:}
\par
(i) {\it Every finite extension of $K$ is a strictly} PQL{\it -field;}
\par
(ii) $P(K) = P$ {\it and $P _{j} (K) = P _{j}$: $j = 0, 1, 2$;}
\par
(iii) {\it For each $p \in (P _{1} \cup P _{2})$, $\gamma _{p}$
equals the degree $[K(\varepsilon _{p}):K]$, where $\varepsilon
_{p}$ is a primitive $p$-th root of unity in $K _{\rm sep}$.}
\par
\medskip
{\it Proof.} Denote by $G _{0}$ and $G _{1}$ the topological group
products $\prod _{p \in (\overline P \setminus P)} \hbox{\Bbb Z}
_{p}$ (i.e. $G _{0} = \{1\}$, in case $P = \overline P$) and
$\prod _{p \in P} \hbox{\Bbb Z} _{p}$, respectively, and fix an
algebraic closure $\overline {\hbox{\Bbb Q}}$ of the field of
rational numbers, as well as a primitive $p$-th root of unity
$\varepsilon _{p} \in \overline {\hbox{\Bbb Q}}$, for each $p \in
\overline P$. Also, let $E _{0}$ be a subfield of $\overline
{\hbox{\Bbb Q}}$, such that $P(E _{0}) = P _{0}$, $P(E _{0}) = P$,
$[E _{0} (\varepsilon _{p}):E _{0}] = \gamma _{p}$: $p \in
(\overline P \setminus P _{0})$, and $G _{E _{0}} \cong G _{1}$
(the existence of $E _{0}$ is guaranteed by [10, Lemma 3.5]).
Suppose further that $\varphi $ is a topological generator of $G
_{E _{0}}$, and for each $p \in (\overline P \setminus P)$,
$\delta _{p}$ is a primitive $\gamma _{p}$-th root of unity in
$\hbox{\Bbb Z} _{p}$, $s _{p}$ and $t _{p}$ are integer numbers,
such that $\varphi (\varepsilon _{p}) = \varepsilon _{p} ^{s
_{p}}$, $t _{p} - \delta _{p} \in p\hbox{\Bbb Z} _{p}$, and $0 \le
s _{p}, t _{p} \le (p - 1)$. Assume also that the roots $\delta
_{p}$ are taken so that $t _{p} = s _{p}$ if and only if $p \in P
_{1}$. Regarding $\hbox{\Bbb Z} _{p}$ as as subgroup of $G _{0}$,
$\forall p \in (\overline P \setminus P)$, consider the
topological semidirect product $G = G _{0} \times G _{E _{0}}$,
defined by the rule $\varphi \lambda _{p}\varphi ^{-1} = \delta
_{p}\lambda _{p}$: $p \in (\overline P \setminus P)$, $\lambda
_{p} \in \hbox{\Bbb Z} _{p}$. It has been proved in [10, Sect. 3]
that there exists a Henselian discrete valued field $(K, v)$, such
that $G _{K}$ is isomorphic to $G$, finite extensions of $K$ are
strictly PQL-fields, $E _{0}$ is a subfield of $\widehat K$, and
$E _{0}$ is algebraically closed in $\widehat K$. In particular,
this implies that $P _{0} (\widehat K) = P _{0}$, $P(\widehat K) =
P$ and $[K(\varepsilon _{p}):K] = \gamma _{p}$: $p \in (\overline
P \setminus P _{0})$. Applying finally (2.2) (ii), one concludes
that $P _{1} (\widehat K) = P _{1}$, and so completes the proof of
Proposition 6.7.
\par
\medskip
{\bf Corollary 6.8.} {\it There exists a set $\{(K _{n}, v _{n}):
n \in \hbox{\Bbb N} \cup {\infty }\}$ of Henselian discrete valued
fields with the following properties:}
\par
(i) {\it The absolute Galois group of a finite extension $R _{n}$
of $K _{n}$ is pronilpotent if and only if $n \in \hbox{\Bbb N}$
and $R _{n}$ contains as a subfield an inertial extension of $K
_{n}$ of degree $n$;}
\par
(ii) {\it Finite extensions of $K _{n}$ are strictly} PQL{\it
-fields, for each $n \in \hbox{\Bbb N}$; they are subject to the
alternative described in Theorem 1.3, provided that $n \ge 2$.}
\par
\medskip
{\it Proof.} This follows at once from Corollary 6.6 and
Proposition 6.7.
\par
\vskip0.6truecm \centerline{\bf 7. Norm groups of formally real
quasilocal fields}
\par
\medskip
In this Section we study the norm groups of formally real
quasilocal fields along the lines drawn in Section 6. It has been
proved in [8, Sect. 3] that a field $K$ is formally real and
quasilocal if and only if it is hereditarily Pythagorean (in the
sense of Becker [4]) with a unique ordering, and cd$_{p}(G _{K})
\le 2$: $p \in P(K(\sqrt{-1})$. Let us note that this occurs if
and only if $2 \not\in P(K(\sqrt{-1}))$ and $G _{K}$ is isomorphic
to the topological semidirect product $G _{K(\sqrt{-1})} \times
\langle \sigma \rangle$, where $G _{K(\sqrt{-1})} \cong \prod _{p
\in P(K(\sqrt{-1}))} \hbox{\Bbb Z} _{p} ^{c(p)}$, $c(p) =$
cd$_{p}(G _{K})$: $p \in P(K(\sqrt{-1}))$, $\sigma ^{2} = 1$, and
$\sigma \tau \sigma ^{-1} = \tau ^{-1}$: $\tau \in G
_{K(\sqrt{-1})}$ (cf. [4, Theorem 1], [5, (3.3)] or [8, (1.2)
and Proposition 3.1]). Note also that if $K$ is formally real and
quasilocal, then $K(\sqrt{-1})$ contains a primitive $m$-th root
of unity, for each $m \in \hbox{\Bbb N}$. Since cyclotomic
extensions are abelian, this can be deduced from the fact (cf. [12,
Lemma 3.6] and [8, Lemma 3.7]) that $K(\sqrt{-1})$ equals
the maximal abelian extension of $K$ in $\overline K$. These results,
combined with the fact that $K$ does not contain a primitive root
of unity of any degree greater than $1$, enable one to obtain
consecutively, and by direct calculations, the following
statements:
\par
\medskip
(7.1) (i) $\sigma (\varepsilon _{n}) = \varepsilon _{n} ^{-1}$,
where $\varepsilon _{n}$ is a primitive $n$-th root of unity in
$K(\sqrt{-1})$, and $n$ is a fixed positive integer not divisible
by $2$;
\par
(ii) With notations being as in (2.2) (ii), $K(\sqrt{-1}) ^{\ast }
= V _{0}$, i.e. $K(\sqrt{-1}) ^{\ast }$ equals the inner group
product $K ^{\ast }K(\sqrt{-1}) ^{\ast p}$, for each prime $p >
2$.
\par
\medskip
These observations will play a crucual role in the proof of the following
statement.
\par
\medskip
{\bf Proposition 7.1.} {\it Let $K$ be a formally real quasilocal
field, $\overline K$ an algebraic closure of $K$, $\overline P$
the set of prime numbers, and $\Pi _{j} = \{p \in \overline P:
$cd$_{p _{j}}(G _{F}) = j\}$, $j = 1, 2$. Assume also that $R$ and
$R _{1}$ are finite extensions of $K$ in $\overline K$, $R _{1} =
R(\sqrt{-1})$, $R _{0}$ is the maximal extension of $K$ in $R$ of
odd degree, and $K _{1} = K(\sqrt{-1})$. Then the following is
true:}
\par
(i) $R$ {\it equals $R _{0}$ or $R _{1}$, and the latter occurs if
and only if $R$ is normal over $K$;}
\par
(ii) $R/K$ {\it possesses an intermediate field $R ^{\prime }$
satisfying the following conditions:}
\par
($\alpha $) {\it the prime divisors of $[R ^{\prime }:K]$ and
$[R:R ^{\prime }]$ lie in $\Pi _{2} \cup \{2\}$ and $\Pi _{1}$,
respectively;}
\par
($\beta $) $N(R/K) = N(R ^{\prime }/K)$ {\it and the group $K
^{\ast }/N(R/K)$ is isomorphic to the direct sum of $G(R ^{\prime
} (\sqrt{-1})/K _{1})$ by a (cyclic) group of order $[R:R _{0}]$.}
\par
\medskip
{\it Proof.} Let $M$ be a finite Galois extension of $K$ including
$R _{1}$. It follows from the structure of $G _{K}$ that $[M:K]$
is even but not divisible by $4$, and $G(M/K)$ possesses a
subgroup $H$ of order $[M:K]/2$. It is also clear that $\varphi
h\varphi ^{-1} = h ^{-1}$: $h \in H$, for every element $\varphi $
of $G(M/K)$ of order $2$. Observing also that $H$ is abelian and
normal in $G(M/K)$, one obtains that subgroups of $G(M/K)$ of odd
orders are included in $H$, and are normal in $G(M/K)$, whereas
subgroups of $G(M/K)$ of even orders equal their normalizers in
$G(M/K)$. Our argument also indicates that if $H _{0}$ is a
subgroup of $G(M/K)$ and $n$ is a positive integer dividing
$[M:K]$ and divisible by the order of $H _{0}$, then there exists
a subgroup $H _{1}$ of $G(M/K)$ of order $n$, such that $H _{0}
\subseteq H _{1}$. These results enable one to deduce Proposition
7.1 (i) from Galois theory, as well as to prove the existence of
an extension $R ^{\prime }$ of $K$ in $R$ satisfying condition
($\alpha $) of Proposition 7.1 (ii). Hence, by Lemma 2.4, we have
$N(R/R ^{\prime }) = R ^{\prime \ast }$, which yields $N(R/K) =
N(R ^{\prime }/K)$. The rest of the proof of the proposition
relies on the fact that the fields $R ^{\prime }$, $R ^{\prime }
\cap R _{0} = R _{0} ^{\prime }$ and $R _{1} ^{\prime } = R _{0}
^{\prime } (\sqrt{-1})$ are related in the same way as $R$, $R
_{0}$ and $R _{1}$. Observing that Br$(K _{1}) _{p}$ is isomorphic
to $\hbox{\Bbb Z} (p ^{\infty })$, for each prime $p$ dividing $[R
_{0} ^{\prime }:K]$ (cf. Proposition 2.5 (iii), [8, (1.2)], and
[30, (11.5)]), one obtains from [11, Theorem 2.1] that $K _{1}
^{\ast }/N(R _{1} ^{\prime }/K _{1})$ is isomorphic to $G(R _{1}
^{\prime }/K _{1})$. Note also that the natural embedding of $R
_{0} ^{\prime }$ into $R _{1}$ induces a group isomorphism $K
^{\ast }/N(R _{0} ^{\prime }/K) \cong K _{1} ^{\ast }/N(R _{1}/K
_{1})$. Indeed, statement (7.1) (ii) implies that $K _{1} ^{\ast }
= K ^{\ast }K _{1} ^{\ast n}$, for every odd integer $n > 2$, so
our assertion reduces to a consequence of Lemma 2.1. Taking
finally into account that $K ^{\ast }/N(R _{1} ^{\prime }/K) \cong
(K ^{\ast }/N(R _{0}/K)) \times (K ^{\ast }/N(K _{1}/K))$, one
completes the proof of Proposition 7.1 (ii).
\par
\medskip
{\bf Remark 7.2.} With assumptions and notations being as in
Proposition 7.1, it is easily seen that a subgroup $H$ of $K
^{\ast }$ is a norm group if and only if the index $\vert K ^{\ast
}:H\vert := n$ is finite and not divisible by any $p \in \Pi
_{1}$. When this occurs, $H$ includes $K ^{\ast e}$, where $e$ is
the exponent of $K ^{\ast }/H$, so it follows from (7.1) (ii),
Proposition 7.1 and the structure of $G _{K}$ (described at the
beginning of the Section) that $e$ is divisible by $n/e$, and $K
^{\ast }/H$ is a direct sum of cyclic groups $C _{e}$ and $C
_{n/e}$ of orders $e$ and $n/e$, respectively.
\par
\medskip
{\bf Corollary 7.3.} {\it For a formally real quasilocal field
$K$, the following assertions hold true:}
\par
(i) {\it In order that $N(R/K) = N(R _{\rm ab}/K)$, for each
finite extension $R$ of $K$, it is necessary and sufficient that
cd$_{p}(G _{K}) \le 1$, for every prime $p \neq 2$;}
\par
(ii) {\it In order that finite extensions of $K$ are strictly}
PQL{\it , it is necessary and sufficient that cd$_{p}(G _{K}) \neq
1$, for each prime $p$; when this occurs, $K$ is either a real
closed field, or $G _{K}$ and finite extensions of $K$ have the
properties required by Theorem 1.3.}
\par
\medskip
{\it Proof.} This follows at once from Proposition 7.1.
\par
\medskip
{\bf Corollary 7.4.} {\it Let $P _{0}$ and $P$ be subsets of the
set $\overline P$ of prime numbers, such that $2 \in P _{0}$ and
$P _{0} \subseteq P$. Then there exist fields $E _{1}$ and $E
_{2}$ with the following properties:}
\par
(i) {\it Finite extensions of $E _{i}$ are strictly} PQL {\it and
$P = \{p \in \overline P$: cd$_{p} (G _{E _{i}}) \neq 0\}$: $i =
1, 2$; moreover, if $P \neq \widetilde P _{j}$, for some index
$j$, where $\widetilde P _{1} = \{2\}$ and $\widetilde P _{2} = P
_{0}$, then $G _{E _{j}}$ is nonnilpotent, and the considered
extensions of $E _{j}$ are subject to the alternative described by
Theorem 1.3;}
\par
(ii) $E _{1}$ {\it is formally real and $E _{2}$ is nonreal with
$P _{0} (E _{2}) = P _{0}$.}
\par
\medskip
{\it Proof.} The existence of $E _{1}$ follows at once from the
classification in [8, Sect. 3] of profinite groups realizable as
absolute Galois groups of formally real quasilocal fields.
Consider now some Henselian discrete valued field $(K, v)$
satisfying the following conditions: (i) char$(\widehat K) = 0$,
$P _{0} (\widehat K) = P _{0}$ and $P _{1} (\widehat K) =
\overline P \setminus P _{0}$; (ii) every finite extension of $K$
is a strictly PQL-field, and the extension of $K$ in $K _{\rm
sep}$, obtained by adjoining a primitive $p$-th root of unity, is
of even degree $\gamma _{p}$, for each $p \in P _{1} (K)$ (the
existence of such $K$ follows from Proposition 6.6). By [7,
Proposition 3.1], $G _{K}$ is a prosolvable group, which means
that it possesses a closed Hall pro-$\Pi $-subgroup $H _{\Pi }$
(uniquely determined, up-to conjugacy in $G _{K}$), for each set
$\Pi $ of prime numbers. In addition, it is easily verified that
the intermediate field $E _{2}$ of $K _{\rm sep}/K$ corresponding
by Galois theory to $H _{P}$ has the properties required by
Corollary 7.4.
\par
\medskip
\vskip0.6truecm
\centerline{\bf Acknowledgements}
\par
\medskip
I would like to thank Professor V. Drensky for communicating to me
some of the references used in the proofs of Theorems 1.1 and 5.1.

\vskip0.9cm \centerline{ REFERENCES} \vglue15pt\baselineskip12.8pt
\def\num#1{\smallskip\item{\hbox to\parindent{\enskip [#1]\hfill}}}
\parindent=1.38cm

\num{1} A. {\pc ALBERT}: {\sl Modern Higher Algebra}. Chicago
Univ. Press, Chicago, Ill., 1937.

\num{2} A. {\pc ALBERT}: {\sl Structure of Algebras}. Amer. Math.
Soc. Colloq. Publ., 24, Providence, RI, 1939.

\num{3} E. {\pc ARTIN}, J. {\pc TATE}: {\sl Class Field Theory}.
Benjamin, New York-Amsterdam, 1968.

\num{4} E. {\pc BECKER}: {\sl Hereditarily Pythagorean fields,
infinite Harrison primes and sums of $2 ^{n}$-th powers}. Bull.
Amer. Math. Soc. {\bf 84} (1978), No 2, 278-280.

\num{5} L. {\pc BR$\ddot o$CKER}: {\sl Characterizations of fans
and hereditarily Pythagorean fields}. Math. Z. {\bf 151} (1976),
149-163.

\num{6} J.W.S. {\pc CASSELS}, A. {\pc FR$\ddot o$HLICH (Eds.)}:
{\sl Algebraic Number Theory}. Academic Press, London-New York,
1967.

\num{7} I.D. {\pc CHIPCHAKOV}: {\sl Henselian valued quasilocal
fields with totally indivisible value groups}. Commun. Algebra
{\bf 27} (1999), 3093-3108.

\num{8} I.D. {\pc CHIPCHAKOV}: {\sl On the Galois cohomological
dimensions of stable fields with Henselian valuations}. Commun.
Algebra {\bf 30} (2002), 1549-1574.

\num{9} I.D. {\pc CHIPCHAKOV}: {\sl Central division algebras of
$p$-primary dimensions and the $p$-component of the Brauer group
of a $p$-quasilocal field}. C.R. Acad. Sci. Bulg. {\bf 55} (2002),
55-60.

\num{10} I.D. {\pc CHIPCHAKOV}: {\sl Henselian discrete valued
fields admitting one-dimensional local class field theory}. Proc.
International Conference on Algebra, Algebraic Geometry and
Applications (V. Brinzanescu, V. Drensky and P. Pragacz, Eds.),
23.9-02.10. 2003, Borovets, Bulgaria, Serdica Math. J. {\bf 30}
(2004), 363-394.

\num{11} I.D. {\pc CHIPCHAKOV}: {\sl One-dimensional abstract
local class field theory}. Preprint.

\num{12} I.D. {\pc CHIPCHAKOV}: {\sl On the residue fields of
Henselian valued stable fields}. Preprint.

\num{13} I.D. {\pc CHIPCHAKOV}: {\sl Algebraic extensions of
global fields admitting one-dimensional local class field theory}.
Preprint.

\num{14} J. {\pc DILLER}, A. {\pc DRESS}: {\sl Zur Galoistheorie
pythagoreischer K$\ddot o$rpern}. Arch. Math. (Basel) {\bf 16}
(1965), 148-152.

\num{15} O. {\pc ENDLER}: {\sl Valuation Theory}. Springer-Verlag,
New York, 1972.

\num{16} B. {\pc FEIN}, M. {\pc SCHACHER}: {\sl Brauer groups of
fields algebraic over $\hbox{\Bbb Q}$}. J. Algebra {\bf 43}
(1976), 328-337.

\num{17} I.B. {\pc FESENKO}: {\sl Local class field theory}.
Perfect residue field case. Izv. Ross. Akad. Nauk, Ser. Mat. {\bf
57} (1993), No 4, 72-91 (Russian: Engl. transl. in Russ. Acad.
Sci., Izv., Math. {\bf 43} (1993), No 4, 72-91).

\num{18} I.B. {\pc FESENKO}, S.V. {\pc VOSTOKOV}: {\sl Local
Fields and Their Extensions}. Transl. Math. Monographs, 121, Am.
Math. Soc., Providence, RI, 2002.

\num{19} L. {\pc FUCHS}: {\sl Infinite Abelian Groups}. Academic
Press, New York-London, 1970.

\num{20} M. {\pc HAZEWINKEL}: {\sl Corps de classes local}.
Appendix to {\it M. Demazure and P. Gabriel}, Groupes Algebriques,
North-Holland, Amsterdam, 1970.

\num{21} K. {\pc IWASAWA}: {\sl Local Class Field Theory}. Iwanami
Shoten, Japan, 1980 (Japanese: Russian transl. by Mir, Moscow,
1983; Engl. transl. in Oxford Mathematical Monographs. Oxford
Univ. Press, New York; Oxford Clarendon Press, VIII, Oxford,
1986).

\num{22} B. {\pc JACOB}, A. {\pc WADSWORTH}: {\sl Division
algebras over Henselian fields}. J. Algebra {\bf 128} (1990),
126-179.

\num{23} M.I. {\pc KARGAPOLOV}, Yu.I. {\pc MERZLYAKOV}: {\sl
Fundamentals of Group Theory, 3rd ed.}, Nauka, Moscow, 1982 (in
Russian).

\num{24} G. {\pc KARPILOVSKY}: {\sl Topics in Field Theory}.
North-Holland Math. Studies, 155, North Holland, Amsterdam etc.,
1989.

\num{25} H. {\pc KOCH}: {\sl Algebraic Number Theory}. Itogi Nauki
i Tekhniki, Ser. Sovrem. Probl. Mat., Fundam. Napravleniya,
Moscow, 1990.

\num{26} S. {\pc LANG}: {\sl Algebra}. Addison-Wesley Publ. Comp.,
Mass., 1965.

\num{27} S. {\pc LANG}: {\sl Algebraic Number Theory}.
Addison-Wesley Publ. Comp., Mass., 1970.

\num{28} P. {\pc MAMMONE}, A. {\pc MERKURJEV}: {\sl On the
corestriction of the $p ^{n}$-symbol}. Commun. Algebra {\bf 76}
(1991), 73-80.

\num{29} A.S. {\pc MERKURJEV}: {\sl Brauer groups of fields}.
Commun. Algebra 11 (1983), No 22, 2611-2624.

\num{30} A.S. {\pc MERKURJEV}, A.A. {\pc SUSLIN}: {\sl
$K$-cohomology of Brauer - Severi varieties and norm residue
homomorphisms}. Izv. Akad. Nauk SSSR {\bf 46} (1982), 1011-1046
(in Russian. Engl. transl. in Math. USSR Izv. {\bf 21} (1983),
307-340).

\num{31} G.A. {\pc MILLER}, H.C. {\pc MORENO}: {\sl Nonabelian
groups in which every subgroup is abelian}. Trans. Amer. Math.
Soc. {\bf 4} (1903), 398-404.

\num{32}  M. {\pc MORIYA}: {\sl Eine notwendige Bedingung f$\ddot
u$r die G$\ddot u$ltigkeit der Klassenk$\ddot o$rper-theorie im
Kleinen}. Math. J. Okayama Univ. {\bf 2} (1952), 13-20.

\num{33} J. {\pc NEUKIRCH}, R. {\pc PERLIS}: {\sl Fields with
local class field theory}. J. Algebra {\bf 42} (1976), 531-536.

\num{34} R. {\pc PIERCE}: {\sl Associative Algebras}.
Springer-Verlag, New York, 1982.

\num{35} V.P. {\pc PLATONOV}, A.S. {\pc RAPINCHUK}: {\sl Algebraic
Groups and Number Theory}. Nauka, Moscow 1991.

\num{36} A. {\pc PRESTEL}, P. {\sl ROQUETTE}: {\sl Formally
$p$-adic Fields}. Lecture Notes in Math. {\bf 1050},
Springer-Verlag, Berlin etc., 1984.

\num{37} L. {\pc REDEI}: {\sl Algebra}, v. 1, Akademiai Kiado,
Budapest, 1967.

\num{38} S.A. {\pc RUSAKOV}: {\sl Analogues to Sylow's theorems on
existence and embeddability of subgroups}. Sibirsk. Mat. Zh. {\bf
4} (1963), No 5, 325-342 (in Russian).

\num{39} O.F.G. {\pc SCHILLING}: {\sl Valuation theory}. Amer.
Math. Soc. Mathematical Surveys, No IV, Maple Press Comp., York,
PA, 1950.

\num{40} O.F.G. {\pc SCHILLING}: {\sl Necessary conditions for
local class field theory}. Math. J. Okayama Univ. {\bf 3} (1953),
5-10.

\num{41} K. {\pc SEKIGUCHI}: {\sl Class field theory of
$p$-extensions over a formal power series field with a
$p$-quasifinite coefficient field}. Tokyo J. Math. {\bf 6} (1983),
167-190.

\num{42} J.-P. {\pc SERRE}: {\sl Cohomologie Galoisienne}. Lecture
Notes in Math. {\bf 5}, Springer-Verlag, Berlin, 1965.

\num{43} J.-P. {\pc SERRE}: {\sl Local Fields}. Graduate Texts in
Mathematics, Springer-Verlag, New York, Heidelberg, Berlin, 1979.

\num{44} I.R. {\pc SHAFAREVICH}: {\sl On $p$-extensions}. Mat.
Sb., n. Ser. {\bf 20(62)} (1947), 351-363 (in Russian: Engl.
transl. in Amer. Math. Soc. Transl., II, Ser. {\bf 4} (1956),
59-72).

\num{45} L.A. {\pc SHEMETKOV}: {\sl Formations of Finite Groups}.
Nauka, Moscow, 1978 (in Russian).

\num{46} J.-P. {\pc TIGNOL}: {\sl On the corestriction of central
simple algebras}. Math. Z. {\bf 194} (1987), 267-274.

\num{47} I.L. {\pc TOMCHIN}, V.I. {\pc YANCHEVSKIJ}: {\sl On
defects of valued division algebras}. Algebra i Analiz {\bf 3}
(1991), No 3, 147-164 (in Russian: Engl. transl. in St. Petersburg
Math. J. {\bf 3} (1992), No 3, 631-646).

\num{48} R. {\pc WARE}: {\sl Galois groups of maximal
$p$-extensions}. Trans. Amer. Math. Soc. {\bf 333} (1992),
721-729.

\num{49} G. {\pc WHAPLES}: {\sl Generalized local class field
theory}. I, Duke Math. J. {\bf 19} (1952), 505-517; II, ibid. {\bf
21} (1954), 247-256; III, ibid., 575-581; IV, ibid. 583-586.

\vskip1cm
\def\pc#1{\eightrm#1\sixrm}
\hfill\vtop{\eightrm\hbox to 5cm{\hfill Ivan {\pc
CHIPCHAKOV}\hfill}
 \hbox to 5cm{\hfill Institute of Mathematics and Informatics\hfill}\vskip-2pt
 \hbox to 5cm{\hfill Bulgarian Academy of Sciences\hfill}
\hbox to 5cm{\hfill Acad. G. Bonchev Str., bl. 8\hfill} \hbox to
5cm{\hfill 1113 {\pc SOFIA,} Bulgaria\hfill}}
\end